\documentclass[12pt]{article}
\usepackage{latexsym, amssymb,amsmath,stmaryrd}
\begin{document}
\baselineskip=16pt

\newcommand{\la}{\langle}
\newcommand{\ra}{\rangle}
\newcommand{\psp}{\vspace{0.4cm}}
\newcommand{\pse}{\vspace{0.2cm}}
\newcommand{\ptl}{\partial}
\newcommand{\dlt}{\delta}
\newcommand{\sgm}{\sigma}
\newcommand{\al}{\alpha}
\newcommand{\be}{\beta}
\newcommand{\G}{\Gamma}
\newcommand{\g}{\gamma}
\newcommand{\gm}{\gamma}
\newcommand{\vs}{\varsigma}
\newcommand{\Lmd}{\Lambda}
\newcommand{\lmd}{\lambda}
\newcommand{\td}{\tilde}
\newcommand{\vf}{\varphi}
\newcommand{\yt}{Y^{\nu}}
\newcommand{\wt}{\mbox{wt}\:}
\newcommand{\der}{\mbox{Der}\:}
\newcommand{\ad}{\mbox{ad}\:}
\newcommand{\stl}{\stackrel}
\newcommand{\ol}{\overline}
\newcommand{\ul}{\underline}
\newcommand{\es}{\epsilon}
\newcommand{\dmd}{\diamond}
\newcommand{\clt}{\clubsuit}
\newcommand{\vt}{\vartheta}
\newcommand{\ves}{\varepsilon}
\newcommand{\dg}{\dagger}
\newcommand{\tr}{\mbox{Tr}\:}
\newcommand{\ga}{{\cal G}({\cal A})}
\newcommand{\hga}{\hat{\cal G}({\cal A})}
\newcommand{\Edo}{\mbox{End}\:}
\newcommand{\for}{\mbox{for}}
\newcommand{\kn}{\mbox{ker}\:}
\newcommand{\Dlt}{\Delta}
\newcommand{\rad}{\mbox{Rad}}
\newcommand{\rta}{\rightarrow}
\newcommand{\mbb}{\mathbb}
\newcommand{\rd}{\mbox{Res}\:}
\newcommand{\stc}{\stackrel{\circ}}
\newcommand{\stdt}{\stackrel{\bullet}}
\newcommand{\lra}{\Longrightarrow}
\newcommand{\f}{\varphi}
\newcommand{\rw}{\rightarrow}
\newcommand{\op}{\oplus}
\newcommand{\co}{\Omega}
\newcommand{\dar}{\Longleftrightarrow}
\newcommand{\sta}{\theta}

\begin{center}{\Large \bf Supersymmetyric Analogues of the Classical}\end{center}
\begin{center}{\Large \bf Theorem on Harmonic Polynomials
} 
\end{center}
\vspace{0.2cm}

\vspace{0.2cm}

\begin{center}{\large Cuiling Luo$^a$ and Xiaoping
Xu$^b$\footnote{Xu's research is supported by NSFC Grant
11171324}}\end{center}

{\small \noindent a. College of Science, Hebei United University,
Tangshan, Hebei 063009, P. R. China.

\noindent b. Corresponding author, Hua Loo-Keng Key Mathematical
Laboratory, Institute of\\ Mathematics, Academy of  Mathematics and
Systems Sciences, Chinese Academy of \\ Sciences, Beijing, 100190,
P. R. China.}

\begin{abstract}{Classical harmonic analysis says that the spaces of
homogeneous harmonic polynomials (solutions of Laplace equation) are
irreducible modules of the corresponding orthogonal Lie group
(algebra) and the whole polynomial algebra is a free module over the
invariant polynomials generated by harmonic polynomials. In this
paper, we first establish two-parameter $\mbb{Z}^2$-graded
supersymmetric oscillator generalizations of the above theorem for
the Lie superalgebra $gl(n|m)$. Then we extend the result to
two-parameter $\mbb{Z}$-graded supersymmetric oscillator
generalizations of the above theorem for the Lie superalgebras
$osp(2n|2m)$ and $osp(2n+1|2m)$.}\end{abstract}

\section{Introduction}

\quad$\;\;$Harmonic polynomials are important objects in analysis,
differential geometry and physics. A fundamental theorem in
classical harmonic analysis says that the spaces of homogeneous
harmonic polynomials (solutions of Laplace equation) are irreducible
modules of the corresponding orthogonal Lie group (algebra) and the
whole polynomial algebra is a free module over the invariant
polynomials generated by harmonic polynomials. Bases of these
irreducible modules can be obtained easily (e.g., cf. [X1]). The
algebraic beauty of the above theorem is that Laplace operator
characterizes the irreducible submodules of the polynomial algebra
and the corresponding quadratic invariant gives a decomposition of
the polynomial algebra into a direct sum of irreducible submodules.

Cao [C] proved that the subspaces of homogeneous polynomial vector
solutions of the $n$-dimensional Navier equations in elasticity are
exactly direct sums of three explicitly given irreducible submodules
when $n\neq 4$ and direct sums of four  explicitly given irreducible
submodules if $n=4$ of the corresponding orthogonal Lie group
(algebra), and the whole polynomial vector space is also a free
module over the invariant polynomials generated these solutions.
This is essentially a vector-function generalization of the
classical theorem on harmonic polynomials.

In [X2], the second author proved that the space of homogeneous
polynomial solutions with degree $m$ for the dual cubic Dickson
invariant differential operator is exactly a direct sum of
$\llbracket m/2 \rrbracket+1$ explicitly determined irreducible
$E_6$-submodules and the whole polynomial algebra is a free module
over the polynomial algebra in the Dickson invariant generated by
these solutions. This gave a cubic $E_6$-generalization of the
classical theorem on harmonic polynomials.

 Lie algebras (Lie groups) serve as the
symmetries in quantum physics (e.g., cf. [FC, L, LF, G]). Their
various representations provide distinct concrete practical physical
models. Many important physical phenomena have been interpreted as
the consequences of symmetry breakings (e.g., cf. [LF]). Harmonic
oscillators are basic objects in quantum mechanics (e.g., cf. [FC,
G]). Oscillator representations of finite-dimensional simple Lie
algebras are the most fundamental ones in quantum physics (e.g., cf.
[DES, FSS]). Howe [Ho] obtained a $\mbb{Z}$-graded
 multiplicity-free oscillator  representation for $sl(n,\mbb{C})$.
In [X1], the second author found the methods of solving flag partial
differential equations for polynomial solutions. Moreover, we [LX]
used a result in [X1] to prove two-parameter $\mbb{Z}^2$-graded
oscillator generalizations of the classical theorem on harmonic
polynomials for the Lie algebra $sl(n,\mbb{C})$ and two-parameter
$\mbb{Z}$-graded oscillator generalizations of the theorem for the
Lie algebra $o(n,\mbb{C})$.

 The aim of this work is to establish two-parameter $\mbb{Z}^2$-graded
supersymmetric oscillator generalizations of the classical theorem
on harmonic polynomials for the Lie superalgebra $gl(n|m)$ and
two-parameter $\mbb{Z}$-graded supersymmetric oscillator
generalizations of the theorem for the Lie superalgebras
$osp(2n|2m)$ and $osp(2n+1|2m)$. Below we give a technical
introduction.

Suppose that $n\geq 3$ is an integer. Denote by $E_{r,s}$ the square
matrix with 1 as its $(r,s)$-entry and 0 as the others. The compact
orthogonal Lie algebra
$$o(n,\mbb{R})=\sum_{1\leq r<s\leq
n}\mbb{R}(E_{r,s}-E_{s,r}).\eqno(1.1)$$ It acts on the polynomial
algebra ${\cal A}=\mbb{R}[x_1,...,x_n]$ by $(E_{r,s}-E_{s,r})|_{\cal
A}=x_r\ptl_{x_s}-x_s\ptl_{x_r}.$
 Recall the Laplace operator
$$\Dlt=\ptl_{x_1}^2+\ptl_{x_2}^2+\cdots+\ptl_{x_n}^2.\eqno(1.2)$$
Moreover, we have the fundamental invariant
$$\eta=x_1^2+x_2^2+\cdots+x_n^2.\eqno(1.3)$$

Denote by ${\cal A}_k$ the subspace of homogeneous polynomials in
${\cal A}$ with degree $k$. Classical theorem on harmonic
polynomials says that the subspace of harmonic polynomials
$${\cal H}_k=\{f\in{\cal A}_k\mid
\Dlt(f)=0\}\eqno(1.4)$$  forms an irreducible $o(n,\mbb{R})$-module
 and
${\cal A}=\bigoplus_{i,k=0}^\infty \eta^i{\cal H}_k$ is a direct sum
of irreducible  $o(n,\mbb{R})$-submodules.  The beauty of the above
theorem is that the invariant differential operator $\Dlt$
characterizes the irreducible submodules and its dual operator
$\eta$ gives the complete reducibility.

Fix two positive integers $m$ and $n$.  Set
$$gl(n|m)_0=\sum_{i,j=1}^n\mbb{C}E_{i,j}+\sum_{r,s=1}^m\mbb{C}E_{n+r,n+s}\eqno(1.5)$$
and
$$gl(n|m)_1=\sum_{i=1}^n\sum_{r=1}^m(\mbb{C}E_{i,n+r}+\mbb{C}E_{n+r,i}).\eqno(1.6)$$
The Lie superalgebra $gl(n|m)=gl(n|m)_0+gl(n|m)_1$ with the
algebraic operation $[\cdot,\cdot]$ defined by
$$[A,B]=AB-(-1)^{i_1i_2}BA\qquad\for\;\;A\in gl(n|m)_{i_1},\;B\in
gl(n|m)_{i_2}.\eqno(1.7)$$

For convenience, we use the notion
$\ol{i,i+j}=\{i,i+1,i+2,...,i+j\}$ for integers $i$ and $j$ with
$i\leq j$. Let ${\cal A}$ be the polynomial algebra in bosonic
variables $\{x_i\mid i\in\ol{1,2n}\}$ and fermionic variables
$\{\sta_j\mid j\in\ol{1,2m}\}$, i.e.,
$$x_rx_s=x_sx_r,\;\sta_p\sta_q=-\sta_q\sta_p,\;\;x_r\sta_p=\sta_px_r\eqno(1.8)
$$ for $r,s\in\ol{1,2n}$ and $p,q\in\ol{1,2m}$. Set
$\Theta=\sum_{p=1}^{2m}\mbb{C}\sta_p.$  Write
$${\cal A}_{(0)}=\sum_{q=0}^m\mbb{C}[x_1,...,x_{2n}]\Theta^{2q},\qquad{\cal
A}_{(1)}=\sum_{q=1}^{m-1}\mbb{C}[x_1,...,x_{2n}]\Theta^{2q+1}.\eqno(1.9)$$
Then ${\cal A}={\cal A}_{(0)}\oplus {\cal A}_{(1)}$ is a
$\mbb{Z}_2$-graded algebra.

For $r\in\ol{1,n}$, the usual differential operator $\ptl_{x_r}$
acts on ${\cal A}$ as a derivation  such that
$\ptl_{x_r}(x_s)=\dlt_{r,s}$ and $\ptl_{x_r}(\sta_p)=0$ for
$s\in\ol{1,2n}$ and $p\in\ol{1,2m}.$  Moreover, for $p\in\ol{1,2m}$,
we define $\ptl_{\sta_p}$ as a linear operator on ${\cal A}$ with
$\ptl_{\sta_p}(x_r)=0$ and $\ptl_{\sta_p}(\sta_q)=\dlt_{p,q}$  for
$r\in\ol{1,n}$ and $q\in\ol{1,2m}$, such that
$$\ptl_{\sta_p}(fg)=\ptl_{\sta_p}(f)g+(-1)^\iota f
\ptl_{\sta_p}(g)\qquad\for\;\;f\in{\cal A}_{(\iota)},\;g\in{\cal
A}.\eqno(1.10)$$

For later notational convenience, we redenote
$$y_i=x_{n+i},\qquad\vt_j=\sta_{m+j}\qquad\for\;\;i\in\ol{1,n},\;j\in\ol{1,m}.\eqno(1.11)$$
 Define a representation of $gl(n|m)$
on ${\cal A}$ determined by
$$E_{i,j}|_{\cal
A}=x_i\ptl_{x_j}-y_j\ptl_{y_i},\qquad E_{i,n+r}|_{\cal
A}=x_i\ptl_{\sta_r}-\vt_r\ptl_{y_i},\eqno(1.12)$$
$$E_{n+r,i}|_{\cal A}=\sta_r\ptl_{x_i}+y_i\ptl_{\vt_r},\qquad
E_{n+r,n+s}|_{\cal
A}=\sta_r\ptl_{\sta_s}-\vt_s\ptl_{\vt_r}\eqno(1.13)$$ for
$i,j\in\ol{1,n}$ and $r,s\in\ol{1,m}$.

Write $\Theta_1=\sum_{r=1}^m\mbb{C}\sta_r$ and
$\Theta_2=\sum_{s=1}^m\mbb{C}\vt_s$. Denote by $\mbb{N}$ the set of
nonnegative integers. For $\ell_1,\ell_2\in\mbb{N}$, we denote
$${\cal A}_{\ell_1,\ell_2}=\mbox{Span}\{x^\al
y^\al\Theta_1^{\ell'_1}\Theta_2^{\ell'_2}\mid\al,\be\in\mbb{N}^n;\ell_1',\ell_2'\in\mbb{N};|\al|+\ell_1'=\ell_1,\;
|\be|+\ell_2'=\ell_2\},\eqno(1.14)$$ where $|\gm|=\sum_{i=1}^n\gm_i$
for $\gm=(\gm_1,...,\gm_n)\in\mbb{N}^n$. Let
$$\Dlt=\sum_{i=1}^n\ptl_{x_i}\ptl_{y_i}+\sum_{r=1}^m\ptl_{\sta_r}\ptl_{\vt_r},\qquad
\eta=\sum_{i=1}^nx_iy_i+\sum_{r=1}^m\sta_r\vt_r.\eqno(1.15)$$
Moreover, we define
$${\cal H}_{\ell_1,\ell_2}=\{f\in{\cal A}_{\ell_1,\ell_2}\mid
\Dlt(f)=0\}.\eqno(1.16)$$ Denote
$k_{\ell,\ell'}=\min\{\ell,\ell'\}.$ The following is our first main
result.\psp

{\bf Theorem 1}. {\it Let $\ell,\ell'\in\mbb{N}$. The space ${\cal
H}_{\ell,\ell'}$ is an irreducible $gl(n|m)$-module  if and only if
$\ell>m+1-n\;\mbox{or}\;\ell'>m+1-n\;\mbox{or}\;\ell+\ell'\leq
m+1-n.$  When $|\ell-\ell'|>m+1-n$ or $\ell+\ell'\leq m+1-n$, ${\cal
A}_{\ell,\ell'}=\bigoplus_{i=0}^{k_{\ell,\ell'}}\eta^i{\cal
H}_{\ell-i,\ell'-i}$ is a decomposition of irreducible
$gl(n|m)$-submodules}.\psp

We remark that if  $\ell,\ell'\leq m+1-n$ and $\ell+\ell'>m+1-n$,
then ${\cal H}_{\ell,\ell'}$ is an indecomposable $gl(n|m)$-module.
In fact, ${\cal H}_{\ell,\ell'}\bigcap\eta{\cal
A}_{\ell-1,\ell'-1}\neq\{0\}$. This also shows that ${\cal
A}_{\ell,\ell'}$ is not completely reducible when $|\ell-\ell'|\leq
m+1-n$ and $\ell+\ell'>m+1-n$.

Fix $1<n_1+1<n_2\leq n$.
 Note
$$[\ptl_{x_r},x_r]=1=[-x_r,\ptl_{x_r}],\;\;[\ptl_{y_s},y_s]=1=[-y_s,\ptl_{y_s}].\eqno(1.17)$$
Changing operators $\ptl_{x_r}\mapsto -x_r,\;
 x_r\mapsto
\ptl_{x_r}$  for $r\in\ol{1,n_1}$ and $\ptl_{y_s}\mapsto -y_s,\;
 y_s\mapsto\ptl_{y_s}$  for $s\in\ol{n_2+1,n}$ in (1.12) and
 (1.13), we get a new representation of $gl(n|m)$ on
  ${\cal A}$ determined by
$$E_{i,j}|_{\cal
A}=E_{i,j}^x-E_{j,i}^y,\qquad E_{n+r,n+s}|_{\cal
A}=\sta_r\ptl_{\sta_s}-\vt_s\ptl_{\vt_r}\eqno(1.18)$$ with
$$E_{i,j}^x|_{\cal A}=\left\{\begin{array}{ll}-x_j\ptl_{x_i}-\delta_{i,j}&\mbox{if}\;
i,j\in\ol{1,n_1};\\ \ptl_{x_i}\ptl_{x_j}&\mbox{if}\;i\in\ol{1,n_1},\;j\in\ol{n_1+1,n};\\
-x_ix_j &\mbox{if}\;i\in\ol{n_1+1,n},\;j\in\ol{1,n_1};\\
x_i\partial_{x_j}&\mbox{if}\;i,j\in\ol{n_1+1,n}
\end{array}\right.\eqno(1.19)$$
and
$$E_{i,j}^y|_{\cal A}=\left\{\begin{array}{ll}y_i\ptl_{y_j}&\mbox{if}\;
i,j\in\ol{1,n_2};\\ -y_iy_j&\mbox{if}\;i\in\ol{1,n_2},\;j\in\ol{n_2+1,n};\\
\ptl_{y_i}\ptl_{y_j} &\mbox{if}\;i\in\ol{n_2+1,n},\;j\in\ol{1,n_2};\\
-y_j\partial_{y_i}-\delta_{i,j}&\mbox{if}\;i,j\in\ol{n_2+1,n};
\end{array}\right.\eqno(1.20)$$
$$ E_{i,n+r}|_{\cal
A}=\left\{\begin{array}{ll}
\ptl_{x_i}\ptl_{\sta_r}-\vt_r\ptl_{y_i}&\mbox{if}\;
i\in\ol{1,n_1};\\ x_i\ptl_{\sta_r}-\vt_r\ptl_{y_i}&\mbox{if}\;
i\in\ol{n_1+1,n_2};\\ x_i\ptl_{\sta_r}+y_i\vt_r&\mbox{if}\;
i\in\ol{n_2+1,n};\end{array}\right.\eqno(1.21)$$
$$ E_{n+r,i}|_{\cal
A}=\left\{\begin{array}{ll}-x_i\sta_r+y_i\ptl_{\vt_r} &\mbox{if}\;
i\in\ol{1,n_1};\\ \sta_r\ptl_{x_i}+y_i\ptl_{\vt_r}&\mbox{if}\; i\in\ol{n_1+1,n_2};\\
\sta_r\ptl_{x_i}+\ptl_{y_i}\ptl_{\vt_r}&\mbox{if}\;
i\in\ol{n_2+1,n}\end{array}\right.\eqno(1.22)$$ for $i,j\in\ol{1,n}$
and $r,s\in\ol{1,m}$.

The related Laplace operator becomes
$$\Dlt=-\sum_{i=1}^{n_1}x_i\ptl_{y_i}+\sum_{r=n_1+1}^{n_2}\ptl_{x_r}\ptl_{y_r}-\sum_{s=n_2+1}^n
y_s\ptl_{x_s}+\sum_{r=1}^m\ptl_{\sta_r}\ptl_{\vt_r}\eqno(1.23)$$ and
its dual
$$\eta=\sum_{i=1}^{n_1}y_i\ptl_{x_i}+\sum_{r=n_1+1}^{n_2}x_ry_r+\sum_{s=n_2+1}^n
x_s\ptl_{y_s}+\sum_{r=1}^m\sta_r\vt_r.\eqno(1.24)$$ Denote
\begin{eqnarray*}\qquad {\cal A}_{\la
\ell_1,\ell_2\ra}&=&\mbox{Span}\{x^\al
y^\be\Theta_1^{\ell_1'}\Theta_2^{\ell_2'}\mid\al,\be\in\mbb{N}^n;\ell_1',\ell_2'\in\mbb{N};\\
& &\sum_{r=n_1+1}^n\al_r-\sum_{i=1}^{n_1}\al_i+\ell_1'=\ell_1;
\sum_{i=1}^{n_2}\be_i-\sum_{r=n_2+1}^n\be_r+\ell_2'=\ell_2\}\hspace{1.8cm}(1.25)\end{eqnarray*}
for $\ell_1,\ell_2\in\mbb{Z}$. Again we set ${\cal
H}_{\la\ell_1,\ell_2\ra}=\{f\in {\cal A}_{\la\ell_1,\ell_2\ra}\mid
\Dlt(f)=0\}. $ The following is our second main result.\psp

{\bf Theorem 2}. {\it Let $\ell,\ell'\in\mbb{Z}$ such that
$\ell'\geq 0$ if $n_2=n$.  The $gl(n|m)$-module ${\cal
H}_{\la\ell,\ell'\ra}$ is irreducible if and only if $\ell+\ell'\leq
n_1+m+1-n_2$ or $\ell\not\in\ol{n_1+1-n,n_1+m+1-n}$ and $n_2=n$.
When $\ell+\ell'\leq n_1+m+1-n_2$, ${\cal
A}_{\la\ell,\ell'\ra}=\bigoplus_{i=0}^\infty\eta^i({\cal
H}_{\la\ell-i,\ell'-i\ra})$ is the decomposition of irreducible
$gl(n|m)$-submodules if  $n_2<n$, and ${\cal
A}_{\la\ell,\ell'\ra}=\bigoplus_{i=0}^{\ell'}\eta^i({\cal
H}_{\la\ell-i,\ell'-i\ra})$ is the decomposition of irreducible
$gl(n|m)$-submodules if $n_2=n$.}\psp

If $\ell+\ell'> n_1+m+1-n_2$  and $\ell\in\ol{n_1+1-n,n_1+m+1-n}$
when $n_2=n$, the $gl(n|m)$-module ${\cal H}_{\la\ell,\ell'\ra}$ is
indecomposable. When $n_2<n$ and $\ell+\ell'> n_1+m+1-n_2$, ${\cal
A}_{\la\ell,\ell'\ra}$ is not completely reducible.

 We use (1.14)
and (1.15) to define
$${\cal A}_k=\bigoplus_{\ell=0}^k{\cal A}_{\ell,k-\ell},\;\;{\cal
H}_k=\{f\in{\cal A}_k\mid\Dlt(f)=0\}\eqno(1.26)$$ for $k\in\mbb{N}$.
Moreover, we use (1.23) and (1.25) to define
$${\cal A}_{\la k\ra}=\bigoplus_{\ell\in\mbb{Z}}{\cal A}_{\la\ell,k-\ell\ra},\;\;{\cal
H}_{\la k\ra}=\{f\in{\cal A}_{\la k\ra}\mid\Dlt(f)=0\}\eqno(1.27)$$
for $k\in\mbb{Z}$. The above representations of $gl(n|m)$ can be
uniquely extended to the representations of the Lie superalgebra
$osp(2n|2m)$. \psp

{\bf Theorem 3}. {\it Suppose $n>1$. For $k\in\mbb{N}$, ${\cal H}_k$
is an irreducible $osp(2n|2m)$-module if and only if $k\leq m+1-n$
or $k>2(m+1-n)$. When $k\leq m+1-n$, ${\cal
A}_k=\bigoplus_{i=0}^{\llbracket k/2 \rrbracket}\eta^i{\cal
H}_{k-2i}$ is a decomposition of irreducible
$osp(2n|2m)$-submodules}.

{\it Let $k\in\mbb{Z}$.  The $osp(2n|2m)$-module ${\cal H}_{\la
k\ra}$ is irreducible if and only if $k\leq n_1+m+1-n_2$. When
$k\leq n_1+m+1-n_2$, ${\cal A}_{\la
k\ra}=\bigoplus_{i=0}^\infty\eta^i({\cal H}_{\la k-2i\ra})$ is the
decomposition of irreducible $osp(2n|2m)$-submodules.}\psp

Let $x_0$ be a bosonic (commuting) variable. Set
$${\cal B}={\cal A}[x_0]=\bigoplus_{k=0}^\infty{\cal B}_k,\qquad{\cal
B}_k=\sum_{i=0}^k{\cal A}_{k-i}x_0^i.\eqno(1.28)$$ Moreover, the
corresponding supersymmetric Laplace operator and invariant become
$$\Dlt'=\ptl_{x_0}^2+2\sum_{i=1}^n\ptl_{x_i}\ptl_{y_i}+2\sum_{r=1}^m\ptl_{\sta_r}\ptl_{\vt_r},\qquad
\eta'=x_0^2+2\sum_{i=1}^nx_iy_i+2\sum_{r=1}^m\sta_r\vt_r.\eqno(1.29)$$
Then the representation of $gl(n|m)$ given in (1.12) and (1.13) can
be uniquely extended to a representation of $osp(2n+1|2m)$ on ${\cal
B}$ such that $\Dlt'$ is an $osp(2n+1|2m)$-invariant operator and
$\eta'$ is an $osp(2n+1|2m)$-invariant. Denote
$${\cal H}'_k=\{f\in{\cal B}_k\mid\Dlt'(f)=0\}\eqno(1.30)$$
for $k\in\mbb{N}$.

Similarly, the representation of $gl(n|m)$ given in (1.18)-(1.22)
can be uniquely extended to a representation of $osp(2n+1|2m)$ on
${\cal B}$ such that the operators
$$\Dlt'=\ptl_{x_0}^2-2(\sum_{i=1}^{n_1}x_i\ptl_{y_i}+\sum_{r=n_1+1}^{n_2}\ptl_{x_r}\ptl_{y_r}-\sum_{s=n_2+1}^n
y_s\ptl_{x_s}+\sum_{r=1}^m\ptl_{\sta_r}\ptl_{\vt_r})\eqno(1.31)$$
and
$$\eta'=x_0^2+2(\sum_{i=1}^{n_1}y_i\ptl_{x_i}+\sum_{r=n_1+1}^{n_2}x_ry_r+\sum_{s=n_2+1}^n
x_s\ptl_{y_s}+\sum_{r=1}^m\sta_r\vt_r)\eqno(1.32)$$ are
$osp(2n+1|2m)$-invariant operators. Set
$${\cal B}_{\la k\ra}=\sum_{i=0}^\infty {\cal A}_{\la
k-i\ra}x_0^i,\qquad {\cal H}'_{\la k\ra}=\{f\in {\cal B}_{\la
k\ra}\mid \Dlt'(f)=0\}.\eqno(1.33)$$\pse

{\bf Theorem 4}.  {\it For any $k\in\mbb{N}$, ${\cal H}_k'$ is an
irreducible $osp(2n+1|2m)$-module. Moreover, ${\cal
B}=\bigoplus_{\ell,k=0}^\infty (\eta')^\ell {\cal H}_k$ is a direct
sum of irreducible $osp(2n+1|2m)$-submodules.}

{\it For any $k\in\mbb{Z}$, ${\cal H}_{\la k\ra}'$ is an irreducible
$osp(2n+1|2m)$-module. Moreover, ${\cal
B}=\bigoplus_{\ell,k=0}^\infty (\eta')^\ell {\cal H}_{\la k\ra}$ is
a direct sum of irreducible $osp(2n+1|2m)$-submodules.}\psp

The first conclusion in Theorem 3 with $n>m+1$ and the first
conclusion in Theorem 4 with $n>m$ were obtained by Zhang [Z].

In Section 2, we give the proof of Theorem 1. We prove Theorem 2 in
Section 3. Section 4 is devoted to the proof of Theorem 3. We show
Theorem 4 in Section 5.

\section{Proof of Theorem 1}

\quad$\;\;$In this section, we want to prove Theorem 1. Recall the
settings in (1.5)-(1.16).

Set
$$\bar{\cal A}=\mbb{C}[x_1,...,x_n,y_1,...,y_n],\qquad \check{\cal
A}=\sum_{i=0}^{2m}\Theta^{2m}.\eqno(2.1)$$ Then $\bar{\cal A}$ and
$\check{\cal A}$ are subalgebras of ${\cal A}$, and ${\cal
A}=\bar{\cal A}\check{\cal A}$. It is straightforward to verify
$$E_{i,j}\Dlt=\Dlt E_{i,j},\qquad E_{i,j}\eta=\eta
E_{i,j}\qquad\for\;\;i,j\in\ol{1,m+n}\eqno(2.2)$$ by (1.12), (1.13)
and (1.15). Indeed, $\eta$ is an invariant, that is,
$E_{i,j}(\eta)=0$ for any $i,j\in\ol{1,m+n}$. Write
$$\bar{\cal G}=\sum_{i,j=1}^n\mbb{C}E_{i,j},\qquad\check{\cal
G}=\sum_{r,s=1}^m\mbb{C}E_{n+r,n+s}.\eqno(2.3)$$ Then they are Lie
subalgebras of $gl(n|m)$. Let
$$\bar H=\sum_{i=1}^n\mbb{C}E_{i,i},\qquad\check H=\sum_{r=1}^m\mbb{C}E_{n+r,n+s},\eqno(2.4)$$
$$\bar{\cal G}_+=\sum_{1\leq i<j\leq n}\mbb{C}E_{i,j},\qquad\check{\cal
G}_+=\sum_{1\leq r<s\leq m}\mbb{C}E_{n+r,n+s}.\eqno(2.5)$$ We take
$\bar H$ as a Cartan subalgebra of $\bar{\cal G}$ and $\bar{\cal
G}_+$ as the subalgebra spanned by positive root vectors in
$\bar{\cal G}$. Similarly, we take $\check H$ as a Cartan subalgebra
of $\check {\cal G}$ and $\check {\cal G}_+$ as the subalgebra
spanned by positive root vectors in $\check {\cal G}$.

Let
$$\bar\Dlt=\sum_{i=1}^n\ptl_{x_i}\ptl_{y_i},\qquad
\bar\eta=\sum_{i=1}^nx_iy_i.\eqno(2.6)$$
 Recall
$$x^\al=x_1^{\al_1}\cdots x_n^{\al_n},\;\;y^\be=y_1^{\be_1}\cdots
y_n^{\be_n}\;\;\for\;\;\al=(\al_1,...,\al_n),\be=(\be_1,...,\be_n)\in\mbb{N}^{\:n}.\eqno(2.7)$$
For $\ell_1,\ell_2\in\mbb{N}$, we denote
$$\bar{\cal A}_{\ell_1,\ell_2}=\mbox{Span}\{x^\al
y^\be\mid\al,\be\in\mbb{N}^{\:n};\sum_{i=1}^n\al_i=\ell_1;\sum_{i=1}^n\be_i=\ell_2\}.\eqno(2.8)$$
Define
$$\bar{\cal H}_{\ell_1,\ell_2}=\{f\in\bar{\cal
A}_{\ell_1,\ell_2}\mid\bar\Dlt(f)=0\}.\eqno(2.9)$$ For
$i\in\ol{1,n}$, we define $\ves_i\in\bar H^\ast$ by
$$\bar\ves_i(E_{j,j})=\dlt_{i,j}\qquad\for\;\;j\in\ol{1,n}.\eqno(2.10)$$
We have the following lemma (e.g., cf. [X1]):\psp

{\bf Lemma 2.1}. {\it Suppose $n>1$. For any $\ell_1,\ell_2$,
$\bar{\cal H}_{\ell_1,\ell_2}$ is a finite-dimensional irreducible
$\bar {\cal G}$-module with  highest-weight vector
$x_1^{\ell_1}y_n^{\ell_2}$ of weight $\ell_1\ves_1+\ell_2\ves_n$.
Moreover,
$$\bar{\cal A}=\bigoplus_{\ell_1,\ell_2,\ell_3=0}^\infty\bar
\eta^{\ell_1}\bar{\cal H}_{\ell_2,\ell_3}\eqno(2.11)$$ is a
decomposition of irreducible $\bar {\cal G}$-submodules.}\psp

When $n=1$, we have
$$\bar{\cal H}_{\ell,0}=\mbb{C}x_1^\ell,\qquad\bar{\cal
H}_{0,\ell}=\mbb{C}y_1^\ell\qquad\mbox{for}\;\;\ell\in\mbb{N}.\eqno(2.12)$$
Moreover,
$$\bar{\cal A}=\bigoplus_{\ell_1,\ell_2=0}^\infty(\bar\eta^{\ell_1}\bar{\cal
H}_{\ell_2,0}\oplus\bar\eta^{\ell_1}\bar{\cal
H}_{0,\ell_2+1}).\eqno(2.13)$$ \psp

Denote
$$\Theta_1=\sum_{i=1}^m\mbb{C}\sta_i,\qquad
\Theta_2=\sum_{i=1}^m\mbb{C}\vt_i.\eqno(2.14)$$ For
$\ell_1,\ell_2\in\ol{1,m}$, we define
$$\check{\cal
A}_{\ell_1,\ell_2}=\Theta_1^{\ell_1}\Theta_2^{\ell_2}.\eqno(2.15)$$
Then $\check{\cal A}_{\ell_1,\ell_2}$ is a finite-dimensional
$\check{\cal G}$-module and
$$\check{\cal
A}=\bigoplus_{\ell_1,\ell_2=0}^m\check{\cal
A}_{\ell_1,\ell_2}.\eqno(2.16)$$ Moreover, we define an ordering:
$$\sta_1\prec\sta_2\prec\cdots\prec\sta_m\prec\vt_m\prec\vt_{m-1}\prec\cdots\prec\vt_1.\eqno(2.17)$$
On the basis
\begin{eqnarray*}\qquad& &\{\sta_{i_1}\cdots\sta_{i_r}\vt_{j_1}\cdots\vt_{j_s}\mid
r,s\in\ol{0,m};\\ & &1\leq i_1<i_2<\cdots<i_r\leq m;m\geq
j_1>j_2>\cdots j_s\geq 1\}\hspace{3.9cm}(2.18)\end{eqnarray*} of
$\check{\cal A}$, we define the partial ordering ``$\prec$"
lexically.

Write
$$\check\Dlt=\sum_{r=1}^m\ptl_{\sta_r}\ptl_{\vt_r},\qquad
\check\eta=\sum_{r=1}^m\sta_r\vt_r.\eqno(2.19)$$ For $r\in\ol{1,m}$,
we define
$$\vec\sta_r=\sta_1\cdots\sta_r,\qquad
\vec\vt_r=\vt_m\cdots\vt_r.\eqno(2.20)$$ For convenience, we let
$$\vec\sta_0=1=\vec\vt_{m+1}.\eqno(2.21)$$
It can easily proved that a minimal term of any singular vector in
$\check{\cal A}_{\ell_1,\ell_2}$ is of the form
$\vec\sta_r\vec\vt_s$ for some $r\in\ol{0,m}$ and $s\in\ol{1,m+1}$
or
$$\vec\sta_r\sta_{r+1}\cdots\sta_{s_1}\vec\vt_{s_2}\vt_{s_1}\cdots\vt_{r+1}\eqno(2.22)$$
for some $0\leq r<s_1<s_2\leq m+1$. A $\check{\cal G}$-{\it singular
vector} $v$ is a nonzero weight vector of $\check{\cal G}$ such that
$\check{\cal G}_+(v)=0$. We count singular vector up to a nonzero
scalar multiple. By comparing minimal terms, we can prove that
$$\{\check\eta^\ell\vec\sta_r\vec\vt_s\mid 0\leq r<s\leq
m+1;\ell\in\ol{0,s-r-1};r+\ell=\ell_1;\ell+m-s+1=\ell_2\}\eqno(2.23)$$
is the set of all $\check{\cal G}$-singular vectors in $\check{\cal
A}_{\ell_1,\ell_2}$. Let $V_{r,s}$ be the finite-dimensional
irreducible $\check{\cal G}$-submodule generated by
$\vec\sta_r\vec\vt_s\in\check{\cal A}_{r,m+1-s}$. By Weyl's Theorem
of complete reducibility (e.g., cf. [Hu]),
$$\check{\cal A}=\bigoplus_{0\leq r<s\leq
m+1}\bigoplus_{\ell=0}^{s-r-1}\check\eta^\ell V_{r,s}\eqno(2.24)$$
is a direct sum of irreducible $\check{\cal G}$-submodules.

Define
$$\check{\cal H}=\{f\in\check{\cal
A}\mid\check\Dlt(f)=0\}.\eqno(2.25)$$ Note
$$E_{n+r,n+s}\check\Dlt=\check\Dlt E_{n+r,n+s},\;\;
E_{n+r,n+s}\check\eta=\check\eta E_{n+r,n+s}\;\;\mbox{on}\;\;{\cal
A} \eqno(2.26)$$ for $r,s\in\ol{1,m}$. Moreover,
$$\check\Dlt\check\eta=\check\eta\check\Dlt-m+\sum_{r=1}^m(\sta_r\ptl_{\sta_r}
+\vt_r\ptl_{\vt_r})\eqno(2.27)$$ by (2.19). Furthermore,
$$\check\Dlt(\vec\sta_r\vec\vt_s)=0\qquad\mbox{if}\;\;r<s.\eqno(2.28)$$
Hence
$$V_{r,s}\subset \check{\cal H}\qquad\for\;\;0\leq r<s\leq
m+1\eqno(2.29)$$ by (2.26). Suppose $0\leq r+1<s\leq m+1$ and
$\ell\in\ol{1,s-r-1}$. For any $f\in V_{r,s}$, we have
$$\check\Dlt(\check\eta^\ell f)=(\sum_{p=0}^{\ell-1}(2p+r+1-s))\check\eta^{\ell-1} f
=\ell(\ell+r-s)\check\eta^{\ell-1} f.\eqno(2.30)$$ Therefore,
$$\check{\cal H}=\bigoplus_{0\leq r<s\leq
m+1} V_{r,s}.\eqno(2.31)$$ In particular,
$$\check{\cal H}_{r,m+1-s}=\{f\in\check{\cal
A}_{r,m+1-s}\mid\check\Dlt(f)=0\}=\check{\cal
A}_{r,m+1-s}\bigcap\check{\cal H}= V_{r,s}\eqno(2.32)$$  for $0\leq
r<s\leq m+1$ and
$$\check{\cal A}_{\ell_1,\ell_2}\bigcap\check{\cal
H}=\{0\}\qquad\mbox{if}\;\;\ell_1+\ell_2\geq m+1.\eqno(2.33)$$
 For
$r\in\ol{1,m}$, we define $\ves_r'\in\check H^\ast$ by
$$\ves_r'(E_{n+s,n+s})=\dlt_{r,s}\qquad\for\;\;s\in\ol{1,m}.\eqno(2.34)$$
Moreover, we treat $\ves'_0=\ves'_{m+1}=0$. Then we have:\psp

{\bf Lemma 2.2}. {\it For $0\leq r<s\leq m+1$, $\check{\cal
H}_{r,m+1-s}$ is a finite-dimensional irreducible $\check{\cal
G}$-module with the highest-weight vector $\vec\sta_r\vec\vt_s$ of
weight $\sum_{p=0}^r\ves_p'-\sum_{q=s}^{m+1}\ves'_q$. Moreover,
$$\check{\cal A}=\bigoplus_{0\leq r<s\leq
m+1}\bigoplus_{\ell=0}^{s-r-1}\check\eta^\ell \check{\cal
H}_{r,m+1-s}.\eqno(2.35)$$}\pse

Recall
$$\bar\Dlt\bar\eta=\bar\eta\bar\Dlt+n+\sum_{i=1}^n(x_i\ptl_{x_i}
+y_i\ptl_{y_i})\eqno(2.36)$$ (e.g., cf. [X1]). Note
$\Dlt=\bar\Dlt+\check\Dlt$. For $\ell_1,\ell_2,\ell_3\in\mbb{N}$,
$0\leq r<s\leq m+1$ and $\ell\in\ol{0,s-r-1}$, $f\in\bar{\cal
H}_{\ell_1,\ell_2}$ and $g\in\check{\cal H}_{r,m+1-s}$, we have
$$\Dlt(\bar\eta^{\ell_3}\check\eta^\ell
fg)=\ell_3(n+\ell_1+\ell_2+\ell_3-1)
\bar\eta^{\ell_3-1}\check\eta^\ell
fg+\ell(\ell+r-s)\bar\eta^{\ell_3}\check\eta^{\ell-1}
fg.\eqno(2.37)$$ Suppose $r+1<s$ and $\ell\in\ol{1,s-r-1}$. If
$$\Dlt(\sum_{p=0}^\ell a_p\bar\eta^p\check\eta^{\ell-p}
fg)=0,\eqno(2.38)$$ then
$$(\ell-p)(p+s-r-\ell)a_p=(p+1)(n+\ell_1+\ell_2+p)a_{p+1}
\qquad\for\;\;p\in\ol{0,\ell-1}.\eqno(2.39)$$ Thus we can take
$a_0=(\ell+1)![\prod_{\iota_2=1}^{\ell+1}
(\iota_2+n+\ell_1+\ell_2-1)]$ and
$$a_{p+1}=[\prod_{\iota_1=0}^p(\ell-\iota_1)(\iota_1+s-r-\ell)][\prod_{\iota_2=p+2}^{\ell+1}
\iota_2(\iota_2+n+\ell_1+\ell_2-1)]\eqno(2.40)$$ for
$p\in\ol{0,\ell-1}$. Denote
\begin{eqnarray*} \Im(\ell_1,\ell_2;r,s,\ell)&=&(\ell+1)![\prod_{\iota_2=1}^{\ell+1}
(\iota_2+n+\ell_1+\ell_2-1)]\check\eta^\ell+\sum_{p=0}^{\ell-1}
[\prod_{\iota_1=0}^p(\ell-\iota_1)(\iota_1+s-r-\ell)]\\
& &\times[\prod_{\iota_2=p+2}^{\ell+1}
\iota_2(\iota_2+n+\ell_1+\ell_2-1)]\bar\eta^{p+1}\check\eta^{\ell-p-1}.\hspace{3.6cm}(2.41)\end{eqnarray*}
 For
convenience, we treat
$$\Im(\ell_1,\ell_2;r,s,0)=n+\ell_1+\ell_2.\eqno(2.42)$$

 Write
$${\cal G}=\bar{\cal G}+\check{\cal G}\cong gl(n,\mbb{C})\oplus
gl(m,\mbb{C}). \eqno(2.43)$$
 Then
$${\cal A}=\bigoplus_{\ell_1,\ell_2,\ell_3=0}^\infty\;\bigoplus_{0\leq r<s\leq
m+1}\bigoplus_{\ell=0}^{s-r-1}(\bar \eta^{\ell_1}\bar{\cal
H}_{\ell_2,\ell_3})(\check\eta^\ell \check{\cal
H}_{r,m+1-s})\eqno(2.44)$$ is a direct sum of irreducible ${\cal
G}$-submodules. By (2.2),
$${\cal H}=\{f\in{\cal A}\mid\Dlt(f)=0\}\eqno(2.45)$$
forms a $gl(n|m)$-submodule, and so it is a ${\cal G}$-submodule.
According to (2.37)-(2.42), $${\cal H}=
\bigoplus_{\ell_2,\ell_3=0}^\infty\;\bigoplus_{0\leq r<s\leq
m+1}\bigoplus_{\ell=0}^{s-r-1} \Im(\ell_2,\ell_3;r,s,\ell)(\bar{\cal
H}_{\ell_2,\ell_3}\check{\cal H}_{r,m+1-s})\eqno(2.46)$$ is a direct
sum of irreducible ${\cal G}$-submodules.

For $\ell,\ell'\in\mbb{N}$, we let
$${\cal
A}_{\ell,\ell'}=\sum_{\ell_1,\ell_2\in\mbb{N},\;\ell_3,\ell_4\in\ol{0,m};\;\ell_1+\ell_3=\ell,\;
\ell_2+\ell_4=\ell'}\bar{\cal A}_{\ell_1,\ell_2}\check{\cal
A}_{\ell_3,\ell_4}\eqno(2.47)$$ and
$${\cal H}_{\ell,\ell'}={\cal A}_{\ell,\ell'}\bigcap{\cal
H}.\eqno(2.48)$$ Then ${\cal A}_{\ell,\ell'}$ and ${\cal
H}_{\ell,\ell'}$ are $gl(n|m)$-submodules. Moreover, $${\cal
H}_{\ell,\ell'}=
\bigoplus_{\ell_1+r+\ell_3=\ell,\;\ell_2+\ell_3+m+1-s=\ell'}
\Im(\ell_1,\ell_2;r,s,\ell_3)\bar{\cal H}_{\ell_1,\ell_2}\check{\cal
H}_{r,m+1-s}\eqno(2.49)$$ is a direct sum of irreducible ${\cal
G}$-submodules. Take the Cartan subalgebra $H=\bar H+\check H$ of
${\cal G}$ (cf. (2.4)) and the subspace
 ${\cal G}_+=\bar{\cal G}_++\check{\cal G}_+$ (cf. (2.5)) spanned by positive
 root vectors in ${\cal G}$.  A ${\cal G}$-{\it singular
vector} $v$ is a nonzero weight vector of ${\cal G}$ such that
${\cal G}_+(v)=0$. We count singular vector up to a nonzero scalar
multiple. Hence we have:\psp

{\bf Lemma 2.3}. {\it The set
\begin{eqnarray*}\qquad&&\{\Im(\ell_1,\ell_2;r,s,\ell_3)(x_1^{\ell_1}y_n^{\ell_2}\vec\sta_r\vec\vt_s)
\mid\ell_1,\ell_2\in\mbb{N};0\leq r<s\leq m+1;\\
& &\ell_3\in\ol{0,s-r-1};\ell_1+r+\ell_3=\ell,\;
\ell_2+\ell_3+m+1-s=\ell'\}\hspace{2.9cm}(2.50)\end{eqnarray*} is
the set of  all the ${\cal G}$-singular vectors in ${\cal
H}_{\ell,\ell'}$, where $\ell_1\ell_2=\ell_1'\ell_2'=0$ if
$n=1$.}\psp

Take $H=\bar H+\check H$ as a Cartan subalgebra of the Lie
superalgebra $gl(n|m)$ and
$$gl(n|m)_+={\cal
G}_++\sum_{r=1}^n\sum_{s=1}^m\mbb{C}E_{r,n+s}\eqno(2.51)$$ as the
subalgebra generated by positive root vectors.  A $gl(n|m)$-{\it
singular vector} $v$ is a nonzero weight vector of $gl(n|m)$ such
that $gl(n|m)_+(v)=0$. We count singular vector up to a nonzero
scalar multiple. Assume that
$x_1^{\ell_1'}y_n^{\ell_2'}\vec\sta_{r'}\vec\vt_{s'}$ is a
$gl(n|m)$-singular vector, where $\ell_1'\ell_2'=0$ when $n=1$. If
$r'\neq 0$, then
\begin{eqnarray*}E_{1,n+r'}(x_1^{\ell_1'}y_n^{\ell_2'}\vec\sta_{r'}\vec\vt_{s'})
&=&
(x_1\ptl_{\sta_{r'}}-\vt_{r'}\ptl_{y_1})(x_1^{\ell_1'}y_n^{\ell_2'}\vec\sta_{r'}\vec\vt_{s'})
\\&=&(-1)^{r'-1}x_1^{\ell_1'+1}y_n^{\ell_2'}\vec\sta_{r'-1}\vec\vt_{s'}-\dlt_{1,n}\ell_2'
x_1^{\ell_1'}y_n^{\ell_2'-1}\vt_{r'}\vec\sta_{r'}\vec\vt_{s'} \neq
0\hspace{1.15cm}(2.52)\end{eqnarray*} by the second equation in
(1.12), which contradicts the definition of singular vector. So
$r'=0$. Suppose $\ell_2'>0$ and $s'>1$. Again the second equation in
(1.12) gives
\begin{eqnarray*}\qquad\quad
E_{n,n+s'-1}(x_1^{\ell_1'}y_n^{\ell_2'}\vec\vt_{s'})&=&
(x_n\ptl_{\sta_{s'-1}}-\vt_{s'-1}\ptl_{y_n})(x_1^{\ell_1'}y_n^{\ell_2'}\vec\vt_{s'})
\\ &=&-\ell_2'
x_1^{\ell_1'}y_n^{\ell_2'-1}\vec\sta_{r'-1}\vec\vt_{s'-1}\neq
0,\hspace{4.5cm}(2.53)\end{eqnarray*} which is absurd. Thus
 $x_1^{\ell_1'}y_n^{\ell_2'}\vec\sta_{r'}\vec\vt_{s'}$ is a
$gl(n|m)$-singular vector if and only if $r'=0$ and
$\ell_2'(s'-1)=0$.

For $\ell_1,\ell_2\in\mbb{N}$, $1\leq r<s-1\leq m$ and
$\ell\in\ol{1,s-r-1}$,
\begin{eqnarray*}&
&E_{1,n+s-1}[\Im(\ell_1,\ell_2;r,s,\ell)]=(x_1\ptl_{\sta_{s-1}}-\vt_{s-1}\ptl_{y_1})[\Im(\ell_1,\ell_2;r,s,\ell)]
\\ &=&x_1\vt_{s-1}\{
(\ell+1)![\prod_{\iota_2=1}^{\ell+1}
(\iota_2+n+\ell_1+\ell_2-1)]\ell \check\eta^{\ell-1}
+\sum_{p=0}^{\ell-2}
[\prod_{\iota_1=0}^p(\ell-\iota_1)(\iota_1+s-r-\ell)]\\
& &\times(\ell-p-1)[\prod_{\iota_2=p+2}^{\ell+1}
\iota_2(\iota_2+n+\ell_1+\ell_2-1)]\bar\eta^{p+1}\check\eta^{\ell-p-2}
\\ & &-\sum_{p=0}^{\ell-1}
(p+1)[\prod_{\iota_1=0}^p(\ell-\iota_1)(\iota_1+s-r-\ell)]
[\prod_{\iota_2=p+2}^{\ell+1}
\iota_2(\iota_2+n+\ell_1+\ell_2-1)]\bar\eta^p\check\eta^{\ell-p-1}\}\\
&=&\ell(\ell+1)(n+\ell+\ell_1+\ell_2+r-s)\Im(\ell_1+1,\ell_2;r,s-1,\ell-1)x_1\vt_{s-1}.\hspace{2.1cm}(2.54)
\end{eqnarray*}
Moreover, (2.41) yields
$$\Im(\ell_1,\ell_2;r,s,\ell)=(\ell+1)![\prod_{\iota_1=0}^\ell(\iota_1+s-r-\ell)]\eta^\ell\;\;\mbox{if}\;\;
n+\ell+\ell_1+\ell_2+r-s=0.\eqno(2.55)$$

Suppose that
$\Im(\ell_1,\ell_2;r,s,\ell_3)(x_1^{\ell_1}y_n^{\ell_2}\vec\sta_r\vec\vt_s)$
is a $gl(n|m)$-singular vector, then
$$n+\ell_3+\ell_1+\ell_2+r-s=0\eqno(2.56)$$
and
$$\Im(\ell_1,\ell_2;r,s,\ell_3)(x_1^{\ell_1}y_n^{\ell_2}\vec\sta_r\vec\vt_s)
=c\eta^{\ell_3}
x_1^{\ell_1}y_n^{\ell_2}\vec\sta_r\vec\vt_s,\eqno(2.57)$$ where
$c=(\ell_3+1)![\prod_{\iota_1=0}^\ell(\iota_1+s-r-\ell_3)]$ by
(2.55). Since
$$E_{i,n+r}(c\eta^{\ell_3}
x_1^{\ell_1}y_n^{\ell_2}\vec\sta_r\vec\vt_s)=c\eta^{\ell_3}
E_{i,n+r}(x_1^{\ell_1}y_n^{\ell_2}\vec\sta_r\vec\vt_s)\eqno(2.58)$$
by (1.12), the arguments in (2.52) and (2.53) show
$r=\ell_2(s-1)=0$. On the other hand, $1\leq 1+r<s$. So $\ell_2=0$.
According to (2.56),
$$\ell_3=s-\ell_1-n.\eqno(2.59)$$
Thus we only get the singular vector
$$\eta^{s-\ell_1-n}
x_1^{\ell_1}\vec\vt_s\in{\cal
H}_{s-n,m+1-n-\ell_1}\;\mbox{with}\;\;s>\ell_1+n.\eqno(2.60)$$

Note $n\leq m$ by $\ell_3>0$. Moreover, $s\leq m+1$ implies
$$s-n,m+1-n-\ell_1\leq
m+1-n.\eqno(2.61)$$ Furthermore,
$$(s-n)+(m+1-n-\ell_1)=m+1-n+\ell_3>m+1-n.\eqno(2.62)$$
This shows that \begin{eqnarray*}\qquad& &{\cal
H}_{\ell,\ell'}\;\mbox{has a unique $gl(n|m)$-singular vector if}\\
& &\ell>m+1-n\;\mbox{or}\;\ell'>m+1-n\;\mbox{or}\;\ell+\ell'\leq
m+1-n.\hspace{3.5cm}(2.63)\end{eqnarray*}

Suppose $n\leq m$ and $\ell,\ell'\in\ol{0,m+1-n}$ such that
$\ell+\ell'>m+1-n$. We take
$$s=n+\ell,\qquad
\ell_1=m+1-n-\ell',\qquad\ell_3=\ell+\ell'+n-m-1.\eqno(2.64)$$ Then
$$\eta^{\ell+\ell'+n-m-1}(x_1^{m+1-n-\ell'}\vec\vt_{n+\ell})\in
{\cal H}_{\ell,\ell'}.\eqno(2.65)$$ Hence $${\cal
H}_{\ell,\ell'}\;\mbox{has exactly two $gl(n|m)$-singular vectors
if}\; \ell+\ell'> m+1-n.\eqno(2.66)$$ In summary, we have: \psp

{\bf Lemma 2.4}. {\it Let $\ell,\ell'\in\mbb{N}$. If
$\ell>m+1-n\;\mbox{or}\;\ell'>m+1-n\;\mbox{or}\;\ell+\ell'\leq
m+1-n,$ the $gl(n|m)$-module ${\cal H}_{\ell,\ell'}$ has a unique
$gl(n|m)$-singular vector. When $\ell,\ell'\leq m+1-n$ and
$\ell+\ell'>m+1-n$, the $gl(n|m)$-module ${\cal H}_{\ell,\ell'}$ has
exactly two $gl(n|m)$-singular vectors.}\psp

Fix $\ell,\ell'\in\mbb{N}$. Let
$$v_{\ell,\ell'}=x_1^\ell y_n^{\ell_1}\vec\vt_s,\qquad
\ell_1+m+1-s=\ell',\;\;\ell_1(s-1)=0,\;\;\dlt_{n,1}\ell\ell_1=0.\eqno(2.67)$$

{\bf Lemma 2.5}. {\it The $gl(n|m)$-module ${\cal H}_{\ell,\ell'}$
is generated by $v_{\ell,\ell'}$.}

{\it Proof.} Let $M$ be the $gl(n|m)$-submodule of ${\cal
H}_{\ell,\ell'}$ generated by $v_{\ell,\ell'}$.  First consider
$n>1$ or $\ell=0$. For any $s'\in\ol{s+1,m+1}$, we have
$$E_{n+s'-1,n}E_{n+s'-2,n}\cdots E_{n+s,n}(v_{\ell,\ell'})=
x_1^\ell y_n^{\ell_1+s'-s}\vec\vt_{s'}\in M\eqno(2.68)$$ by (1.13).
In other words, \begin{eqnarray*}\qquad\quad & &x_1^\ell
y_n^{\ell_2}\vec\vt_{s'}\in M\;\;\mbox{for
any}\;\ell_2\in\mbb{N}\;\mbox{and}\;s'\in\ol{1,m+1}\\ & &\mbox{such
that}\;\;\ell_2+m+1-s'=\ell'\;\mbox{and}\;\dlt_{n,1}\ell\ell_2=0.
\hspace{4.6cm}(2.69)\end{eqnarray*} If $\ell >0$, then for any
$1\leq r\leq \mbox{min}\{\ell,s'-1\}$, we have
$$E_{n+1,1}E_{n+2,1}\cdots E_{n+r,1}(x_1^\ell
y_n^{\ell_2}\vec\vt_{s'})=[\prod_{p=0}^{r-1}(\ell-p)]x_1^{\ell-r}y_n^{\ell_2}\vec\sta_r\vec\vt_{s'}\in
M\eqno(2.70)$$ by (1.13) again. Thus we have showed that
$$x_1^{\ell_3}y_n^{\ell_2}\vec\sta_r\vec\vt_{s'}\in
M\;\;\mbox{whenever}\;\;r+\ell_3=\ell,\;\ell_2+m+1-s'=\ell',\;\;\dlt_{n,1}\ell_2\ell_3=0\eqno(2.71)$$
for $\ell_2,\ell_3\in\mbb{N}$ and $0\leq r<s'\leq m+1$. Recall the
Lie algebra ${\cal G}$ defined in (2.43). As ${\cal G}$-modules,
$$\sum_{r+\ell_3=\ell,\;\ell_2+m+1-s'=\ell'}\check{\cal
H}_{\ell_3,\ell_2}\bar{\cal H}_{r,m+1-s'}\subset M.\eqno(2.72)$$

For $i\in\mbb{N}+1$, we define
$${\cal H}^{(i)}_{\ell,\ell'}=\{f\in{\cal H}_{\ell,\ell'}\mid
\bar\Dlt^i(f)=0.\}\eqno(2.73)$$ Then
$${\cal H}^{(1)}_{\ell,\ell'}=\sum_{r+\ell_3=\ell,\;\ell_2+m+1-s'=\ell'}\check{\cal
H}_{\ell_3,\ell_2}\bar{\cal H}_{r,m+1-s'}\subset M.\eqno(2.74)$$
Denote
$$k_{\ell,\ell'}=\min\{\ell,\ell'\}.\eqno(2.75)$$
Then
$${\cal H}_{\ell,\ell'}^{(k_{\ell,\ell'}+1)}={\cal
H}_{\ell,\ell'}\eqno(2.76)$$ by (2.49). Let
$\ell_1,\ell_2,\in\mbb{N}$, $0\leq r'+1<s'\leq m+1$ and
$\ell_3\in\ol{1,s'-r'-1}$ such that $\ell_1+\ell_3+r'=\ell$ and
$\ell_2+\ell_3+m+1-s'=\ell'$. Then
$$\Im(\ell_1,\ell_2;r',s',\ell_3)(x_1^{\ell_1}y_n^{\ell_2}\vec\sta_{r'}\vec\vt_{s'})
\in {\cal H}_{\ell,\ell'}^{(\ell_3+1)}\eqno(2.77)$$ and
$$\bar\Dlt^{\ell_3}\Im(\ell_1,\ell_2;r',s',\ell_3)
(x_1^{\ell_1}y_n^{\ell_2}\vec\sta_{r'}\vec\vt_{s'})=
c(\ell_1,\ell_2;r',s',\ell_3)x_1^{\ell_1}y_n^{\ell_2}\vec\sta_{r'}\vec\vt_{s'}\eqno(2.78)$$
with
$$c(\ell_1,\ell_2;r',s',\ell_3)=
\ell_3(n+\ell_1+\ell_3+\ell_3-1)(n+\ell_1+\ell_2+\ell_3)(\ell_3+1)!
[\prod_{p=1}^{\ell_3}(s'-r'-p)]\eqno(2.79)$$ by (2.36) and (2.41).
On the other hand,
$$x_1^{\ell_1+\ell_3}y_n^{\ell_2}\vec\sta_{r'}\vec\vt_{s'-\ell_3}\in
M\;\;\mbox{with}\;\;\dlt_{n,1}\ell_2=0\eqno(2.80)$$ by (2.71) and
\begin{eqnarray*}\qquad & &M\ni f=E_{n+s'-1,1}E_{n+s'-2,1}\cdots
E_{n+s'-\ell_3,1}(x_1^{\ell_1+\ell_3}y_n^{\ell_2}\vec\sta_{r'}\vec\vt_{s'-\ell_3})
\\ &=&(-1)^{r'\ell_3}x_1^{\ell_1+\ell_3}y_1^{\ell_3}
y_n^{\ell_2}\vec\sta_{r'}\vec\vt_{s'}+y_n^{\ell_2}\sum_{i=0}^{\ell_3-1}\zeta_iy_1^i
\hspace{6.3cm}(2.81)\end{eqnarray*} with
$\zeta_0,...,\zeta_{\ell_3-1}\in\mbb{C}[x_1]\check{\cal A}$ (cf.
(2.3)). Moreover,
$$\bar\Dlt^{\ell_3}[(\ell_3!)^2{\ell_1+\ell_3\choose\ell_3}\Im(\ell_1,\ell_2;r',s',\ell_3)
(x_1^{\ell_1}y_n^{\ell_2}\vec\sta_{r'}\vec\vt_{s'})-(-1)^{r'\ell_3}
c(\ell_1,\ell_2;r',s',\ell_3)f]=0.\eqno(2.82)$$ Hence
$$\Im(\ell_1,\ell_2;r',s',\ell_3)(\bar{\cal
H}_{\ell_1,\ell_2}\check{\cal H}_{r',m+1-s'})\subset {\cal
H}^{(\ell_3)}_{\ell,\ell'}+M.\eqno(2.83)$$ By (2.49) and induction
on $i$, we have
$${\cal H}^{(i)}_{\ell,\ell'}\subset M\qquad\mbox{for
any}\;i\in\mbb{N}+1.\eqno(2.84)$$ According to (2.86), ${\cal
H}_{\ell,\ell'}=M$. So ${\cal H}_{\ell,\ell'}$ is generated by
$v_{\ell,\ell'}.\qquad\Box$\psp

{\it Proof of Theorem 1}\psp

 According to (2.66), a necessary condition
for ${\cal H}_{\ell,\ell'}$ to be an irreducible $gl(n|m)$-module is
$\ell>m+1-n\;\mbox{or}\;\ell'>m+1-n\;\mbox{or}\;\ell+\ell'\leq
m+1-n.$ To prove the sufficiency, we suppose that
$\ell>m+1-n\;\mbox{or}\;\ell'>m+1-n\;\mbox{or}\;\ell+\ell'\leq
m+1-n.$ Let $V$ be a nonzero $gl(n|m)$-submodule of ${\cal
H}_{\ell,\ell'}$. According to Lemma 2.4, ${\cal H}_{\ell,\ell'}$
has a unique singular vector $v_{\ell,\ell'}$ (cf. (2.67)). Since
$V$ is finite-dimensional, it contains a singular vector. So
$v_{\ell,\ell'}\in V$. Lemma 2.5 says $V={\cal H}_{\ell,\ell'}$.
Hence ${\cal H}_{\ell,\ell'}$ is an irreducible $gl(n|m)$-module.

Let $\ell_1,\ell_2\in\mbb{N}$ and $0\leq r<s\leq m+1$ such that
$n+\ell_1+\ell_2+r-s\geq 0$ and $\dlt_{n,1}\ell_1\ell_2=0$.
 For any
$\ell_4\in\mbb{N}+1$ and $\ell_3\in\ol{0,s-r-1}$, $$\Dlt^{\ell_4+1}
(\eta^{\ell_4}\Im(\ell_1,\ell_2;r,s,\ell_3)x_1^{\ell_1}y_n^{\ell_2}\vec\sta_r\vec\vt_s)=0\eqno(2.85)$$
and
\begin{eqnarray*} & &\Dlt^{\ell_4}
(\eta^{\ell_4}\Im(\ell_1,\ell_2;r,s,\ell_3)x_1^{\ell_1}y_n^{\ell_2}\vec\sta_r\vec\vt_s)
\\ &=&\ell_4![\prod_{i=1}^{\ell_4}(n+i+\ell_1+2\ell_3+r-s)]
\Im(\ell_1,\ell_2;r,s,\ell_3)(x_1^{\ell_1}y_n^{\ell_2}\vec\sta_r\vec\vt_s)\neq
0\hspace{2.6cm}(2.86)\end{eqnarray*} by (2.27) and (2.36). Thus the
set
$$\{\eta^{\ell_4}\Im(\ell_1,\ell_2;r,s,\ell_3)x_1^{\ell_1}y_n^{\ell_2}\vec\sta_r\vec\vt_s\mid\ell_4\in\mbb{N},
\;\ell_3\in\ol{0,s-r-1}\}\eqno(2.87)$$ is linearly independent.

Note that
$$\check\eta^{s-r}x_1^{\ell_1}y_n^{\ell_2}\vec\sta_r\vec\vt_s=0\eqno(2.88)$$
by (2.19) and (2.20). So for any $k\in\mbb{N}$,
\begin{eqnarray*}& &\mbox{Span}\{\eta^{\ell_4}\Im(\ell_1,\ell_2;r,s,\ell_3)x_1^{\ell_1}y_n^{\ell_2}
\vec\sta_r\vec\vt_s\mid\ell_4\in\mbb{N},\;\ell_3\in\ol{0,s-r-1};\ell_3+\ell_4=k\}
\\ & \subset& \mbox{Span}\{
\bar\eta^{\ell_5}\check\eta^{\ell_6}x_1^{\ell_1}y_n^{\ell_2}
\vec\sta_r\vec\vt_s\mid\ell_5\in\mbb{N};\ell_6\in\ol{0,s-r-1};\ell_5+\ell_6=k\}.
\hspace{2.9cm}(2.89)\end{eqnarray*} But the linear independency of
(2.87) implies that the above subspaces have the same dimension.
Thus
\begin{eqnarray*}& &\mbox{Span}\{\eta^{\ell_4}\Im(\ell_1,\ell_2;r,s,\ell_3)x_1^{\ell_1}y_n^{\ell_2}
\vec\sta_r\vec\vt_s\mid\ell_4\in\mbb{N},\;\ell_3\in\ol{0,s-r-1}\}
\\ & =& \mbox{Span}\{
\bar\eta^{\ell_5}\check\eta^{\ell_6}x_1^{\ell_1}y_n^{\ell_2}
\vec\sta_r\vec\vt_s\mid\ell_5\in\mbb{N};\ell_6\in\ol{0,s-r-1}\}.
\hspace{5cm}(2.90)\end{eqnarray*} Therefore, as ${\cal G}$-modules,
$$\bigoplus_{\ell_3=0}^{s-r-1}\sum_{\ell_4=0}^\infty
\eta^{\ell_4}\Im(\ell_1,\ell_2;r,s,\ell_3)\bar{\cal
H}_{\ell_1,\ell_2}\check{\cal H}_{r,m+1-s}
=\bigoplus_{\ell_5=0}^{s-r-1}\sum_{\ell_6=0}^\infty
\bar\eta^{\ell_5}\check\eta^{\ell_6}\bar{\cal
H}_{\ell_1,\ell_2}\check{\cal H}_{r,m+1-s}.\eqno(2.91)$$

Assume that $|\ell-\ell'|>m+1-n$ or $\ell+\ell'\leq m+1-n$.
According to (2.44) and (2.90), the ${\cal G}$-module
\begin{eqnarray*} &{\cal A}_{\ell,\ell'}=
\mbox{Span}\{\bar\eta^{\ell_5}\check\eta^{\ell_6}\bar{\cal
H}_{\ell_1,\ell_2}\check{\cal
H}_{r,m+1-s}\mid\ell_1,\ell_2,\ell_6\in\mbb{N};0\leq r<s\leq m+1;
\ell_5\in\ol{0,s-r-1};\\
&\dlt_{n,1}\ell_1\ell_2=0;\ell_1+\ell_5+\ell_6+r=\ell,\;\ell_2+\ell_5+\ell_6+m+1-s=\ell'\}
\\&=\mbox{Span}\{\eta^{\ell_4}\Im(\ell_1,\ell_2;r,s,\ell_3)\bar{\cal
H}_{\ell_1,\ell_2}\check{\cal
H}_{r,m+1-s}\mid\ell_1,\ell_2,\ell_6\in\mbb{N};\dlt_{n,1}\ell_1\ell_2=0; 0\leq r<s\leq m+1;\\
&
\ell_5\in\ol{0,s-r-1};\ell_1+\ell_3+\ell_4+r=\ell,\;\ell_2+\ell_3+\ell_4+m+1-s=\ell'\}.
\hspace{2cm}(2.92)\end{eqnarray*}

According to  (2.49), (2.91) and (2.92),
$${\cal A}_{\ell,\ell'}=\bigoplus_{i=0}^{k_{\ell,\ell'}}\eta^i{\cal
H}_{\ell-i,\ell'-i}\eqno(2.94)$$ (cf. (2.75)). Let
$i\in\ol{0,k_{\ell,\ell'}}$. If $|\ell-\ell'|>m+1-n$, then
$$\ell-i\geq |\ell-\ell'|>m+1-n\;\;\mbox{or}\;\;\ell'-i\geq
|\ell-\ell'|>m+1-n.\eqno(2.95)$$ When $\ell+\ell'\leq m+1-n$,
$$(\ell-i)+(\ell'-i)=\ell+\ell'-2i
\leq m+1-n.\eqno(2.96)$$ Thus all ${\cal H}_{\ell-i,\ell'-i}$ are
irreducible $gl(n|m)$-submodules. Hence (2.94) is a direct sum of
irreducible $gl(n|m)$-submodules.

This completes the proof of Theorem 1. $\qquad\Box$ \psp

The following result will be used to obtain explicit bases for
modules.\psp

{\bf Lemma 2.6 (Xu [X1])}. {\it Let ${\cal B}$ be a commutative
associative algebra and let ${\cal A}$ be a free ${\cal B}$-module
 generated  by a filtrated subspace $V=\bigcup_{r=0}^\infty V_r$
(i.e., $V_r\subset V_{r+1}$). Let $T_1$ be a linear operator on
${\cal A}$ with a right inverse $T_1^-$ such that
$$T_1({\cal B}),\;T_1^-({\cal B})\subset{\cal B},\qquad
T_1(\eta_1\eta_2)=T_1(\eta_1)\eta_2,\qquad
T_1^-(\eta_1\eta_2)=T_1^-(\eta_1)\eta_2 \eqno(2.97)$$ for $\eta_1
\in {\cal B},\;\eta_2\in V$, and let $T_2$ be a linear operator on
${\cal A}$ such that
$$ T_2(V_{r+1})\subset {\cal B}V_r,\;\;
T_2(f\zeta)=fT_2(\zeta) \qquad\for\;\; r\in\mbb{N},\;\;f\in{\cal
B},\;\zeta\in{\cal A}.\eqno(2.98)$$ Then we
have \begin{eqnarray*}\hspace{1cm}&&\{g\in{\cal A}\mid (T_1+T_2)(g)=0\}\\
& =&\mbox{Span}\{ \sum_{i=0}^\infty(-T_1^-T_2)^i(hg)\mid g\in
V,\;h\in {\cal B};\;T_1(h)=0\}. \hspace{3.7cm}(2.99)\end{eqnarray*}
} \psp

Set
$$\es_i=(0,...,0,\stl{i}{1},0,...,0)\in \mbb{N}^{\:n}.\eqno(2.100)$$
 For each
$i\in\ol{1,n}$, we define the linear operator $\int_{(x_i)}$ on
${\cal A}$ by:
$$\int_{(x_i)}(x^\al)=\frac{x^{\al+\es_i}}{\al_i+1}\;\;\for\;\;\al\in
\mbb{N}^{\:n}.\eqno(2.101)$$ Furthermore, we let
$$\int_{(x_i)}^{(0)}=1,\qquad\int_{(x_i)}^{(m)}=\stl{m}{\overbrace{\int_{(x_i)}\cdots\int_{(x_i)}}}
\qquad\for\; \;0<m\in\mbb{Z}\eqno(2.102)$$ and denote
$$\ptl^{\al}=\ptl_{x_1}^{\al_1}\ptl_{x_2}^{\al_2}\cdots
\ptl_{x_n}^{\al_n},\;\;
\int^{(\al)}=\int_{(x_1)}^{(\al_1)}\int_{(x_2)}^{(\al_2)}\cdots
\int_{(x_n)}^{(\al_n)}\qquad\for\;\;\al\in
\mbb{N}^{\:n}.\eqno(2.103)$$ Obviously, $\int^{(\al)}$ is a right
inverse of $\ptl^\al$ for $\al\in \mbb{N}^{\:n}.$ We remark that
$\int^{(\al)}\ptl^\al\neq 1$ if $\al\neq 0$ due to $\ptl^\al(1)=0$.
In this paper,  our $T_1$'s are of the type $\ptl^\al$ and the right
inverse $T_1^-=\int^{(\al)}$.

Denote $\G_0=\emptyset$ and
$$\G_\ell=\{\vec j=(j_1,j_2,..,j_\ell)\mid 1\leq j_1<j_2<\cdots<
j_\ell\leq m\}\qquad\for\;\;\ell\in\ol{1,m}.\eqno(2.104)$$ Moreover,
we set
$$\sta_\emptyset=\vt_\emptyset=1,\;\;\sta_{\vec
j}=\sta_{j_1}\sta_{j_2}\cdots\sta_{j_\ell},\;\;\vt_{\vec
j}=\vt_{j_1}\vt_{j_2}\cdots\vt_{j_\ell}.\eqno(2.105)$$ Then the set
\begin{eqnarray*}& &\{\sum_{i=0}^\infty\frac{(-1)^ix_1^iy_1^i}{\prod_{r=1}^i(\al_1+i)(\be_1+i)}
(\sum_{s=2}^n\ptl_{x_2}\ptl_{y_2}+\sum_{r=1}^m\ptl_{\sta_r}\ptl_{\vt_r})^i(x^\al
y^\be\sta_{\vec j}\vt_{\vec k})\mid\al,\be\in\mbb{N}^n;\\ & &\vec
j\in\G_{\ell_1};\vec
k\in\G_{\ell_2};\al_1\be_1=0;\ell_1,\ell_2\in\ol{0,m};|\al|+\ell_1=\ell;|\be|+\ell_2=\ell'\}\hspace{2.4cm}(2.106)\end{eqnarray*}
forms a basis of ${\cal H}_{\ell,\ell'}$ by Lemma 2.6 with
$T_1=\ptl_{x_1}\ptl_{y_1},\;T_2=\Dlt-T_1$ and
$T_1^-=\int_{(x_1)}\int_{(y_1)}$.\psp

{\bf Remark 2.7}. If  $\ell,\ell'\leq m+1-n$ and $\ell+\ell'>m+1-n$,
then ${\cal H}_{\ell,\ell'}$ is an indecomposable $gl(n|m)$-module
by (2.66) and Lemma 2.5, and ${\cal H}_{\ell,\ell'}\bigcap\eta{\cal
A}_{\ell-1,\ell'-1}\neq\{0\}$. This also shows that ${\cal
A}_{\ell,\ell'}$ is not completely reducible when $|\ell-\ell'|\leq
m+1-n$ and $\ell+\ell'>m+1-n$.

\section{Proof of Theorem 2}

\quad$\;\;$In this section, we want to prove Theorem 2. Recall the
settings in (1.5)-(1.11) and (1.18)-(1.25).

The Laplace operator in (2.6) changes to
$$\bar\Dlt=-\sum_{i=1}^{n_1}x_i\ptl_{y_i}+\sum_{r=n_1+1}^{n_2}\ptl_{x_r}\ptl_{y_r}-\sum_{s=n_2+1}^n
y_s\ptl_{x_s}\eqno(3.1)$$ and its dual changes to
$$\bar\eta=\sum_{i=1}^{n_1}y_i\ptl_{x_i}+\sum_{r=n_1+1}^{n_2}x_ry_r+\sum_{s=n_2+1}^n
x_s\ptl_{y_s}.\eqno(3.2)$$ We take (2.19) and then supersymmetric
Laplace operator in (1.23) and its dual in (1.24) can be written as:
$$\Dlt=\bar\Dlt+\check\Dlt,\qquad
\eta=\bar\eta+\check\eta.\eqno(3.3)$$ Then with respect to the
representation in (1.18)-(1.22), we have
$$E_{i,j}\Dlt=\Dlt E_{i,j},\qquad E_{i,j}\eta=\eta
E_{i,j}\qquad\for\;\;i,j\in\ol{1,m+n}.\eqno(3.4)$$ Moreover, we take
the settings in (2.1), (2.3)-(2.5) and  (2.7).

Denote
$$\bar{\cal A}_{\la \ell_1,\ell_2\ra}=\mbox{Span}\{x^\al
y^\be\mid\al,\be\in\mbb{N}\:^n;\sum_{r=n_1+1}^n\al_r-\sum_{i=1}^{n_1}\al_i=\ell_1;
\sum_{i=1}^{n_2}\be_i-\sum_{r=n_2+1}^n\be_r=\ell_2\}\eqno(3.5)$$ for
$\ell_1,\ell_2\in\mbb{Z}$. Then $\bar{\cal A}_{\la
\ell_1,\ell_2\ra}$ forms a $\bar{\cal G}$-submodule. Moreover, for
$\ell,\ell'\in\mbb{Z}$, we let
$${\cal
A}_{\la\ell,\ell'\ra}=\sum_{\ell_1,\ell_2\in\mbb{Z},\;\ell_3,\ell_4\in\ol{0,m};\;\ell_1+\ell_3=\ell,\;
\ell_2+\ell_4=\ell'}\bar{\cal A}_{\la\ell_1,\ell_2\ra}\check{\cal
A}_{\ell_3,\ell_4}.\eqno(3.6)$$ It can be verified that ${\cal
A}_{\la\ell,\ell'\ra}$ forms a $gl(n|m)$-submodule. Define
$${\cal H}=\{f\in {\cal A}\mid \Dlt(f)=0\},\;\;{\cal H}_{\la\ell_1,\ell_2\ra}={\cal
H}\bigcap{\cal A}_{\la\ell_1,\ell_2\ra}.\eqno(3.7)$$ By (3.4),
${\cal H}_{\la\ell,\ell'\ra}$ forms a $gl(n|m)$-submodule of ${\cal
A}_{\la\ell,\ell'\ra}$. According to [LX], we have:\psp

{\bf Lemma 3.1}. {\it The  nonzero vectors in
$$\{\mbb{C}[\bar\eta](x_i^{m_1}y_j^{m_2})\mid
m_1,m_2\in\mbb{N};i=n_1,n_1+1;j=n_2,n_2+1-\dlt_{n_2,n}\}\eqno(3.8)$$
are all the singular vectors of $\bar{\cal G}$ (cf. (2.3)-(2.5)) in
$\bar{\cal A}$ (cf. (2.1))}.\psp

Recall ${\cal G}=\bar{\cal G}+\check{\cal G}$ (cf. (2.3)). Take the
Cartan subalgebra $H=\bar H+\check H$ of ${\cal G}$ (cf. (2.4)) and
the subspace
 ${\cal G}_+=\bar{\cal G}_++\check{\cal G}_+$ (cf. (2.5)) spanned by positive
 root vectors in ${\cal G}$. Then
the  nonzero vectors in
\begin{eqnarray*}\hspace{2cm}& &\{\mbb{C}[\bar\eta,\check\eta](x_i^{m_1}y_j^{m_2}\vec\sta_r\vec\vt_s)\mid
m_1,m_2\in\mbb{N};i=n_1,n_1+1;\\ & &j=n_2,n_2+1-\dlt_{n_2,n};0\leq
r<s\leq m+1\}\hspace{4.7cm}(3.9)\end{eqnarray*} are all the ${\cal
G}$-singular vectors   in ${\cal A}$.
 Choose $H$ as a Cartan
subalgebra of the Lie superalgebra $gl(n|m)$ and $gl(n|m)_+={\cal
G}_++\sum_{r=1}^n\sum_{s=1}^m\mbb{C}E_{r,n+s}$ as the subalgebra
generated by positive root vectors.

Fix $\ell,\ell'\in\mbb{Z}$. Then a $gl(n|m)$-singular vector in
${\cal A}_{\la\ell,\ell'\ra}$ must be of the form
$$f=\sum_{p=0}^{\min\{s-r-1,\ell_1\}}b_p\bar\eta^{\ell_1-p}\check\eta^p
(x_i^{m_1}y_j^{m_2}\vec\sta_r\vec\vt_s),\eqno(3.10)$$ where
$\ell_1,m_1,m_2\in\mbb{N},\;\;0\leq r<s\leq m+1,\;b_p\in\mbb{C}$ and
$$(i,j)\in\{(n_1,n_2),
(n_1,n_2+1-\dlt_{n_2,n}),(n_1+1,n_2),(n_1+1,n_2+1-\dlt_{n_2,n})\}.\eqno(3.11)$$

Suppose that $\ell_1=0$ or $s-r-1=0$. Then we can assume
$f=x_i^{m_1}y_j^{m_2}\vec\sta_r\vec\vt_s$. If $r\neq 0$, then (1.21)
implies
$$E_{n_1+1,n+r}(f)=(-1)^{r-1}x_{n_1+1}x_i^{m_1}y_j^{m_2}\vec\sta_{r-1}\vec\vt_s\neq 0,\eqno(3.12)$$
which is absurd. So $r=0$. Assume that $m_2>0,\;s>1$ and $j=n_2$.
According to (1.21),
$$E_{n_2,n+1}(f)=-m_2x_i^{m_1}y_j^{m_2-1}\vt_1\vec\vt_s\neq
0,\eqno(3.13)$$ which leads a contradiction. Hence $m_2(s-1)=0$ if
$j=n_2$. If $n_2<n$, then we have
$$E_{n_2+1,n+1}(f)=x_i^{m_1}y_j^{m_2}y_{n_2+1}\vt_1\vec\vt_s\neq
0\eqno(3.14)$$ by (1.21), which is absurd. Thus $s=1$. In summary,
$$f=x_i^{m_1}y_j^{m_2}\vec\vt_s\qquad\mbox{with}\;s=1\;\mbox{or}\;n_2=n\;\mbox{and}\;m_2=0.\eqno(3.15)$$

Consider the case $\ell_1>0$ and $s-r-1>0$. By (1.21) and the fact
$\check\eta^{s-r-1}\vt_{r+1}\vec\sta_r\vec\vt_s=0$, we have
\begin{eqnarray*}\qquad 0&=&E_{n_1+1,n+r+1}(f)=(x_{n_1+1}\ptl_{\sta_{r+1}}-\vt_{r+1}\ptl_{y_{n_1+1}})(f)
\\
&=&[\sum_{p=1}^{\min\{s-r-1,\ell_1\}}pb_p\bar\eta^{\ell_1-p}\check\eta^{p-1}-
\sum_{p=0}^{\min\{s-r-1,\ell_1\}-1}(\ell_1-p)b_p\bar\eta^{\ell_1-p-1}\check\eta^p
(x_i^{m_1}y_j^{m_2}\vec\sta_r\vec\vt_s)]\\& &\times
x_{n_1+1}\vt_{r+1}x_i^{m_1}y_j^{m_2}\vec\sta_r\vec\vt_s,
\hspace{8.6cm}(3.16)\end{eqnarray*} which implies
$$f=b_0(\bar\eta+\check\eta)^{\ell_1}(x_i^{m_1}y_j^{m_2}\vec\sta_r\vec\vt_s)=b_0\eta^{\ell_1}(x_i^{m_1}y_j^{m_2}\vec\sta_r\vec\vt_s).
\eqno(3.17)$$ The arguments in the previous paragraph and (3.4)
give:\psp

{\bf Lemma 3.2}. {\it Any $gl(n|m)$-singular vector in ${\cal
A}_{\la\ell,\ell'\ra}$ must be of the form:
$$\eta^{\ell_1}(x_i^{m_1}y_j^{m_2}\vec\vt_s)\;\;\mbox{with}\;\ell_1,m_1,m_2\in\mbb{N},\;s\in\ol{1,m+1}\;\mbox{and (3.10)}
\eqno(3.18)$$ such that $s=1$ or $n_2=n$ and $m_2=0$.}\psp

We define
$$\flat=\sum_{r=n_1+1}^nx_r\ptl_{x_r}-\sum_{i=1}^{n_1}x_i\ptl_{x_i},\;\;
\flat'=\sum_{i=1}^{n_2}y_i\ptl{y_i}-\sum_{r=n_2+1}^ny_r\ptl{y_r}.\eqno(3.19)$$
Then
$$\bar{\cal A}_{\la\ell,\ell'\ra}=\{f\in\bar{\cal A}\mid
\flat(f)=\ell f;\flat'(f)=\ell' f\}.\eqno(3.20)$$
 We calculate
$$\bar\Dlt\bar\eta=\bar\eta\bar\Dlt+n_2-n_1+\flat+\flat'.\eqno(3.21)$$
For $\ell_1\in\mbb{N}+1$ and $f\in{\cal H}_{\la\ell,\ell'\ra}$ (cf.
(3.7)), (2.27) and (3.21) imply
$$\Dlt\eta^{\ell_1}(f)=\ell_1(n_2-n_1-m+\ell+\ell'+\ell_1-1)\eta^{\ell_1-1}(f).\eqno(3.22)$$
Thus
$$\Dlt\eta^{\ell_1}(f)=0\dar \ell+\ell'\leq
n_1+m-n_2\;\;\mbox{and}\;\;\ell_1=n_1+m-n_2-\ell-\ell'+1.\eqno(3.23)$$
If the condition holds, then
$$\eta^{\ell_1}(f)\in{\cal H}_{\la
n_1+m-n_2-\ell'+1,n_1+m-n_2-\ell+1\ra}.\eqno(3.24)$$ Moreover,
$$(n_1+m-n_2-\ell'+1)+(n_1+m-n_2-\ell+1)\geq
n_1+m-n_2+2.\eqno(3.25)$$

Observe that
$$x_{n_1}^{m_1}y_{n_2}^{m_2}\vec\vt_s\in{\cal
H}_{\la-m_1,m+1+m_2-s\ra},\qquad
x_{n_1}^{m_1}y_{n_2+1}^{m_2}\vec\vt_1\in{\cal
H}_{\la-m_1,m-m_2\ra},\eqno(3.27)$$
$$x_{n_1+1}^{m_1}y_{n_2}^{m_2}\vec\vt_s\in{\cal
H}_{\la m_1,m+1+m_2-s\ra},\qquad
x_{n_1+1}^{m_1}y_{n_2+1}^{m_2}\vec\vt_1\in{\cal H}_{\la
m_1,m-m_2\ra}.\eqno(3.28)$$ By Lemma 3.2,
$$\mbox{any nonzero}\;{\cal H}_{\la\ell,\ell'\ra}\;\;\mbox{contains a singular
vector of the form}\;x_i^{m_1}y_j^{m_2}\vec\vt_s,\eqno(3.29)$$ where
$s=1$ or $n_2=n$ and $m_2=0$.

Now we consider $f=x_i^{m_1}y_j^{m_2}\vec\vt_s$ with
$m_1,m_2\in\mbb{N},\;s\in\ol{1,m+1}$ and (3.11) such that $s=1$ or
$n_2=n$ and $m_2=0$. Assume $\Dlt\eta^{\ell_1}(f)=0$ for some
$\ell_1\in\mbb{N}+1$.\psp

{\it Case 1}. $(i,j)=(n_1,n_2)$.\psp

In this subcase,
 $\ell=-m_1$ and
$\ell'=m_2+m+1-s$ by (3.27). Thus $m_2-m_1+1-s\leq n_1-n_2$ and
$\ell_1= n_1+m_1+s-n_2-m_2$ by (3.23). So
$$\eta^{\ell_1}(f)\in{\cal H}_{\la
n_1+s-n_2-m_2,n_1+m_1-n_2+m+1\ra}.\eqno(3.30)$$

{\it Case 2}. $(i,j)=(n_1,n_2+1)$.\psp

In this subcase $s=1$, $\ell=-m_1$ and $\ell'=m-m_2$ by (3.27). Thus
$m_1+m_2\geq n_2-n_1$ and $ \ell_1=n_1+m_1+m_2-n_2+1$ by (3.23).
Hence
$$\eta^{\ell_1}(f)\in{\cal H}_{\la
n_1+m_2-n_2+1,n_1+m_1-n_2+m+1\ra}.\eqno(3.31)$$

{\it Case 3}. $(i,j)=(n_1+1,n_2)$.\psp

In this subcase,
 $\ell=m_1$ and
$\ell'=m_2+m+1-s$ by (3.28). Thus $m_2+m_1+1-s\leq n_1-n_2$ and
$\ell_1= n_1-m_1+s-n_2-m_2$ by (3.23). So
$$\eta^{\ell_1}(f)\in{\cal H}_{\la
n_1+s-n_2-m_2,n_1-m_1-n_2+m+1\ra}.\eqno(3.32)$$

{\it Case 4}. $(i,j)=(n_1+1,n_2+1)$.\psp

In this subcase $s=1$, $\ell=m_1$ and $\ell'=m-m_2$ by (3.28). Thus
$m_1-m_2\leq n_1-n_2$ and $ \ell_1=n_1-m_1+m_2-n_2+1$ by (3.23).
Hence
$$\eta^{\ell_1}(f)\in{\cal H}_{\la
n_1+m_2-n_2+1,n_1-m_1-n_2+m+1\ra}.\eqno(3.33)$$

Thus we obtain:\psp

{\bf Lemma 3.3}. {\it A nonzero $gl(n|m)$-module ${\cal
H}_{\la\ell,\ell'\ra}$ has a unique singular vector if and only if
$\ell+\ell'\leq n_1+m+1-n_2$ or $\ell\not\in\ol{n_1+1-n,n_1+m+1-n}$
and $n_2=n$. If the condition holds, the unique singular vector is
of the form $x_i^{m_1}y_j^{m_2}\vec\vt_s$ with (3.11), where $s=1$
or $n_2=n$ and $m_2=0$.}\psp

Fix ${\cal H}_{\la\ell,\ell'\ra}\neq \{0\}$. Assume
$$v_{\ell,\ell'}=x_i^{m_1}y_j^{m_2}\vec\vt_s\in {\cal
H}_{\la\ell,\ell'\ra}\eqno(3.34)$$ for some $(i,j)$ in (3.11),
$m_1,m_2\in\mbb{N}$ and $s\in\ol{1,m+1}$ such that $s=1$ or $n_2=n$
and $m_2=0$.\psp

{\bf Lemma 3.4}. {\it As a $gl(n|m)$-module, ${\cal
H}_{\la\ell,\ell'\ra}$ is generated by $v_{\ell,\ell'}$.}

{\it Proof}. Denote
$$\td\G=\{\td\al=(\al_{n_1+1},...,\al_{n_2})\in\mbb{N}^{n_2-n_1}\},\;\;|\td\al|=\sum_{i=1}^{n_2-n_1}\al_{n_1+i}.
\eqno(3.35)$$ Set
$$\td{\cal A}=\check{\cal
A}[x_{n_1+1},...,x_{n_2},y_{n_1+1},...,y_{n_2}],\eqno(3.36)$$
$$\td{\cal H}_{\la\ell_1,\ell_2\ra}=\td{\cal A}\bigcap{\cal
H}_{\la \ell_1,\ell_2\ra},\eqno(3.37)$$
$$\td\Dlt=-\sum_{i=1}^{n_1}x_i\ptl_{y_i}+\sum_{r=n_1+2}^{n_2}\ptl_{x_r}\ptl_{y_r}-\sum_{s=n_2+1}^n
y_s\ptl_{x_s}+\check\Dlt.\eqno(3.38)$$ Then
$\Dlt=\ptl_{x_{n_1+1}}\ptl_{y_{n_1+1}}+\td\Dlt$. For any
$k_1,k_2\in\mbb{N}$, we define the operator
$$T_{k_1,k_2}=\sum_{i=0}^\infty
\frac{(-1)^ix_{n_1+1}^{k_1+i}y_{n_1+1}^{k_2+i}}{\prod_{r=1}^i(k_1+r)(k_2+r)}\td\Dlt^i.\eqno(3.39)$$
Then Lemma 2.6 yields
\begin{eqnarray*}\qquad{\cal
H}_{\la\ell,\ell'\ra}&=&\mbox{Span}\{T_{\al_{n_1+1},\be_{n_1+1}}([\prod_{i\neq
n_1+1}x_i^{\al_i}y_i^{\be_i}])\sta_{\vec j}\vt_{\vec
k})\mid\al,\be\in\mbb{N}^n;\\ & &\qquad\vec j\in\G_{\ell_1};\vec
k\in\G_{\ell_2};\ell_1,\ell_2\in\ol{0,m};\al_{n_1+1}\be_{n_1+1}=0;\\
& &\qquad\sum_{s=n_1+1}^n\al_s-\sum_{r=1}^{n_1}\al_r+\ell_1=\ell;
\sum_{r=1}^{n_2}\be_r-\sum_{s=n_2+1}^n\be_s+\ell_2=\ell'\},\hspace{1.1cm}(3.40)\end{eqnarray*}
\begin{eqnarray*}\td{\cal
H}_{\la\ell,\ell'\ra}\!\!\!&=&\!\!\!\mbox{Span}\{T_{\al_{n_1+1},\be_{n_1+1}}([\prod_{i=n_1+2}^{n_2}
x_i^{\al_i}y_i^{\be_i}])\sta_{\vec j}\vt_{\vec
k})\mid\td\al,\td\be\in\td\G;\vec j\in\G_{\ell_1};\vec
k\in\G_{\ell_2};\\ & &\qquad
\al_{n_1+1}\be_{n_1+1}=0;\ell_1,\ell_2\in\ol{0,m};|\td\al|+\ell_1=\ell;
|\td\be|+\ell_2=\ell'\}.\hspace{2.3cm}(3.41)\end{eqnarray*}

Write
$${\cal G}_1=\sum_{i,j=1}^{n_1}\mbb{C}E_{i,j},\qquad {\cal
G}_2=\sum_{r,s=n_2+1}^n\mbb{C}E_{r,s},\eqno(3.42)$$
$${\cal G}_3=\sum_{n_1+1\neq i,j\in\ol{1,n_2}}\mbb{C}E_{i,j},\qquad {\cal
G}_4=\sum_{r,s=n_1+2}^n\mbb{C}E_{r,s}.\eqno(3.43)$$ Then
$$\xi\td\Dlt=\td\Dlt\xi\qquad\for\;\;\xi\in {\cal
G}_i,\;i\in\ol{1,4}.\eqno(3.44)$$
 Denote by $V$
the $gl(n|m)$-submodule of ${\cal H}_{\la\ell,\ell'\ra}$ generated
by $v_{\ell,\ell'}$. By (1.18)-(1.20)
$$-E_{n_1+1,n_1}|_{\cal
A}=x_{n_1}x_{n_1+1}+y_{n_1}\ptl_{y_{n_1+1}},\;\;E_{n_2+1,n_2}|_{\cal
A}=x_{n_2+1}\ptl_{x_{n_2}}+y_{n_2}y_{n_2+1}.\eqno(3.45)$$ According
to (1.21) and (1.22),
$$E_{n_2,n+r}=-x_{n_2}\ptl_{\sta_r}+\vt_r\ptl_{y_{n_2}},\;\;E_{n+r,n_2}=\sta_r\ptl_{x_{n_2}}+y_{n_2}\ptl_{\vt_r}
\qquad\for\;\;r\in\ol{1,m}.\eqno(3.46)$$ Repeatedly applying the
operators in (3.45) and (3.46) to $v_{\ell,\ell'}$, we obtain
$$x_{n_1}^{p_1}x_{n_1+1}^{p_2}y_{n_2}^{p_3}y_{n_2+1}^{p_4}\vec\vt_{s'}\in
V\eqno(3.47)$$ for $p_i\in\mbb{N}$ and $s'\in\ol{1,m+1}$ such that
$$p_2-p_1=\ell,\;p_3-p_4+m+1-s'=\ell',\;p_3(s'-1)=0.\eqno(3.48)$$
Lemma 2.5, (2.67),  and (3.37)  tell us that
$$x_{n_1}^{p_1}y_{n_2+1}^{p_4}\td{\cal H}_{\la
p_2,p_3+m+1-s'\ra}\subset V.\eqno(3.49)$$

Let $U({\cal G}_i)$ be the universal enveloping of the Lie algebra
${\cal G}_i$. Applying $U({\cal G}_1)$ and $U({\cal G}_2)$ to
(3.49), we get
$$[\prod_{\iota_1=1}^{n_1}x_{\iota_1}^{\al_{\iota_1}}][\prod_{\iota_2=n_2+1}^n y_{\iota_2}^{\be_{\iota_2}}]
\td{\cal H}_{\la p_2,p_3+m+1-s'\ra}\subset V\eqno(3.50)$$ for any
$(\al_1,...,\al_{n_1})\in\mbb{N}^{n_1}$ and
$(\be_{n_2+1},...,\be_n)\in\mbb{N}^{n-n_2}$ such that
$\sum_{\iota_1=1}^{n_1}\al_{\iota_1}=p_1$ and
$\sum_{\iota_2=n_2+1}^n \be_{\iota_2}=p_4$. Applying $U({\cal G}_3)$
to (3.50), we obtain
$$T_{\al_{n_1+1},\be_{n_1+1}}([\prod_{i\neq
n_1+1}x_i^{\al_i}y_i^{\be_i}])\sta_{\vec j}\vt_{\vec k})\in
V\eqno(3.51)$$ by (3.41) and (3.44), where $\al,\be,\vec j,\vec k$
are as those in (3.40) and $\al_{n_2+1}=\cdots=\al_n=0$. Finally, we
get (3.51) with any  $\al,\be,\vec j,\vec k$ in (3.40) by (3.44) and
applying $U({\cal G}_4)$. According to (3.40), $V={\cal
H}_{\la\ell,\ell'\ra}.\qquad\Box$\psp

{\it Proof of Theorem 2}\psp

 Suppose $\ell+\ell'\leq n_1+m+1-n_2$ or
$\ell\not\in\ol{n_1+1-n,n_1+m+1-n}$ and $n_2=n$. Let $V$ be a
nonzero submodule of ${\cal H}_{\la\ell,\ell'\ra}$. According to
Lemma 3.3, the vector $v_{\ell,\ell'}$ in (3.34) is the unique
singular vector of ${\cal H}_{\la\ell,\ell'\ra}$. Since $gl(n|m)_+$
in (2.51) is locally nilpotent by (1.18)-(1.22), $V$ contains a
singular vector. So $v_{\ell,\ell'}\in V$. By Lemma 3.4, $V={\cal
H}_{\la\ell,\ell'\ra}$, that is, ${\cal H}_{\la\ell,\ell'\ra}$ is
irreducible. The necessity also follows from Lemma 3.3.

Assume $\ell+\ell'\leq n_1+m+1-n_2$. Since $\Dlt$ is locally
nilpotent by (3.1) and (3.3), for any $0\neq u\in {\cal
A}_{\la\ell,\ell'\ra}$, there exists an element
$\kappa(u)\in\mbb{N}$ such that
$$\Dlt^{\kappa(u)}(u)\neq
0\;\;\mbox{and}\;\;\Dlt^{\kappa(u)+1}(u)=0.\eqno(3.52)$$ Set
$$\Psi=\left\{\begin{array}{ll}\sum_{i=0}^\infty\eta^i({\cal
H}_{\la\ell-i,\ell'-i\ra})&\mbox{if}\;n_2<n,\\
\sum_{i=0}^{\ell'}\eta^i({\cal
H}_{\la\ell-i,\ell'-i\ra})&\mbox{if}\;n_2=n.\end{array}\right.\eqno(3.53)$$
Given $0\neq u\in {\cal A}_{\la\ell,\ell'\ra}$, $\kappa(u)=1$
implies $u\in {\cal H}_{\la\ell,\ell'\ra}\subset\Psi$. Suppose that
$u\in \Psi$ whenever $\kappa(u)<r$ for some positive integer $r$.
Assume $\kappa(u)=r$. First
$$v=\Dlt^r(u)\in {\cal
H}_{\la\ell-r,\ell'-r\ra}\subset\Psi.\eqno(3.54)$$ Note
$$\Dlt^r[\eta^r(v)]=r![\prod_{i=1}^r(n_2-n_1-m+\ell+\ell'-r-i)]v\eqno(3.55)$$
by (3.22). Thus we  have either
$$u=\frac{1}{r![\prod_{i=1}^r(n_2-n_1-m+\ell+\ell'-r-i)]}\eta^r(v)\in\Psi\eqno(3.56)$$
or
$$\kappa\left(u-\frac{1}{r![\prod_{i=1}^r(n_2-n_1-m+\ell+\ell'-r-i)]}\eta^r(v)\right)<r.\eqno(3.57)$$
By induction,
$$u-\frac{1}{r![\prod_{i=1}^r(n_2-n_1-m+\ell+\ell'-r-i)]}\eta^r(v)\in\Psi,\eqno(3.58)$$
which implies $u\in\Psi$. Therefore, we have $\Psi={\cal
A}_{\la\ell,\ell'\ra}$. Since all $\eta^i({\cal
H}_{\la\ell-i,\ell'-i\ra})$ have distinct highest weights, the sums
in (3.53) are direct sums.

This completes the proof of Theorem 2. $\qquad\Box$\psp

{\bf Remark 3.6}. If $\ell+\ell'> n_1+m+1-n_2$  and $\ell\leq
n_1+m+1-n_2$ when $n_2=n$, the $gl(n|m)$-module ${\cal
H}_{\la\ell,\ell'\ra}$ is indecomposable. When $\ell+\ell'>
n_1+m+1-n_2$,  ${\cal A}_{\la\ell,\ell'\ra}$ is not completely
reducible.\psp

Recall the notations in (2.104) and (2.105). The set
\begin{eqnarray*}\hspace{1cm}& &\{\sum_{i=0}^\infty\frac{(-1)^ix_{n_1+1}^iy_{n_1+1}^i}
{\prod_{r=1}^i(\al_{n_1+1}+r)(\be_{n_1+1}+r)}
(\Dlt-\ptl_{x_{n_1+1}}\ptl_{y_{n_1+1}})^i(x^\al y^\be\sta_{\vec
j}\vt_{\vec k})\\ & &\qquad \mid\al,\be\in\mbb{N}^n;\vec
j\in\G_{\ell_1};\vec
k\in\G_{\ell_2};\al_{n_1+1}\be_{n_1+1}=0;\ell_1,\ell_2\in\ol{0,m};\\
& &\qquad \sum_{r=n_1+1}^n\al_r-\sum_{i=1}^{n_1}\al_i+\ell_1=\ell;
\sum_{i=1}^{n_2}\be_i-\sum_{r=n_2+1}^n\be_r+\ell_2=\ell'\}\hspace{2.5cm}(3.59)\end{eqnarray*}
 forms a basis of ${\cal H}_{\la\ell,\ell'\ra}$ by Lemma 2.6.

\section{Proof of Theorem 3}

\quad$\;\;$In this section, we want to prove Theorem 3.

 Set
$${\cal K}_0=\sum_{i,j=1}^n\mbb{C}(E_{i,j}-E_{n+j,n+i})+\sum_{r,s=1}^m\mbb{C}
(E_{2n+r,2n+s}-E_{2n+m+s,2n+m+r}),\eqno(4.1)$$
$${\cal
K}_1=\sum_{i=1}^n\sum_{r=1}^m[\mbb{C}(E_{i,2n+r}-E_{2n+m+r,n+i})+\mbb{C}
(E_{2n+r,i}+E_{n+i,2n+m+r})].\eqno(4.2)$$ Then ${\cal K}={\cal
K}_0+{\cal K}_1$ forms a Lie sub-superalgebra of $gl(2n|2m)$
isomorphic to $gl(n|m)$. Let \begin{eqnarray*} \qquad
osp(2n|2m)_0&=&{\cal K}_0+\sum_{1\leq i<j\leq
n}[\mbb{C}(E_{i,n+j}-E_{j,n+i})+\mbb{C}(E_{n+i,j}-E_{n+j,i})]
\\ & &+\sum_{1\leq r\leq s\leq m}[\mbb{C} (E_{2n+r,2n+m+s}+E_{2n+s,2n+m+r})
\\ & &\quad+\mbb{C}
(E_{2n+m+r,2n+s}+E_{2n+m+s,2n+r})],\hspace{4.3cm}(4.3)\end{eqnarray*}
$$osp(2n|2m)_1={\cal
K}_1+\sum_{i=1}^n\sum_{r=1}^m[\mbb{C}(E_{i,2n+m+r}+E_{2n+r,n+i})+\mbb{C}
(E_{n+i,2n+r}-E_{2n+m+r,i})].\eqno(4.4)$$ The space
$osp(2n|2m)=osp(2n|2m)_0+osp(2n|2m)_1$ forms a simple Lie
sub-superalgebra of of $gl(2n|2m)$. Moreover, its Lie subalgebra
$$osp(2n|2m)_0\cong o(2n,\mbb{C})\oplus sp(2n,\mbb{C}).\eqno(4.5)$$

Take settings in (1.8)-(1.11).  Define a representation of
$osp(2n|2m)$ on ${\cal A}$ determined by
$$(E_{i,j}-E_{n+j,n+i})|_{\cal
A}=x_i\ptl_{x_j}-y_j\ptl_{y_i},\;\;
(E_{2n+r,2n+s}-E_{2n+m+s,2n+m+r})|_{\cal
A}=\sta_r\ptl_{\sta_s}-\vt_s\ptl_{\vt_r},\eqno(4.6)$$
$$(E_{i,2n+r}-E_{2n+m+r,n+i})|_{\cal
A}=x_i\ptl_{\sta_r}-\vt_r\ptl_{y_i},\;\;
(E_{2n+r,i}+E_{n+i,2n+m+r})|_{\cal
A}=\sta_r\ptl_{x_i}+y_i\ptl_{\vt_r},\eqno(4.7)$$
$$(E_{i,n+j}-E_{j,n+i})|_{\cal A}=x_i\ptl_{y_j}-x_j\ptl_{y_i},\;(E_{2n+m+r,2n+s}+E_{2n+m+s,2n+r})]|_{\cal
A}=\vt_r\ptl_{\sta_s}+\vt_s\ptl_{\sta_r},\eqno(4.8)$$
$$(E_{n+i,j}-E_{n+j,i})|_{\cal
A}=y_i\ptl_{x_j}-y_j\ptl_{x_i},\;\;(E_{2n+r,2n+m+s}+E_{2n+s,2n+m+r})|_{\cal
A}=\sta_r\ptl_{\vt_s}+\sta_s\ptl_{\vt_r},\eqno(4.9)$$
$$(E_{i,2n+m+r}+E_{2n+r,n+i})|_{\cal
A}=x_i\ptl_{\vt_r}+\sta_r\ptl_{y_i},\;
(E_{n+i,2n+r}-E_{2n+m+r,i})|_{\cal
A}=y_i\ptl_{\sta_r}-\vt_r\ptl_{x_i}\eqno(4.10)$$
 for $i,j\in\ol{1,n}$ and $r,s\in\ol{1,m}$.

Recall that we write $\Theta_1=\sum_{r=1}^m\mbb{C}\sta_r$ and
$\Theta_2=\sum_{s=1}^m\mbb{C}\vt_s$. For $k\in\mbb{N}$, we denote
$${\cal A}_k=\mbox{Span}\{x^\al
y^\al\Theta_1^{\ell'_1}\Theta_2^{\ell'_2}\mid\al,\be\in\mbb{N}^n;\ell_1',\ell_2'\in\mbb{N};|\al|+\ell_1'+
|\be|+\ell_2'=k\}.\eqno(4.11)$$ Again we take
$$\Dlt=\sum_{i=1}^n\ptl_{x_i}\ptl_{y_i}+\sum_{r=1}^m\ptl_{\sta_r}\ptl_{\vt_r},\qquad
\eta=\sum_{i=1}^nx_iy_i+\sum_{r=1}^m\sta_r\vt_r.\eqno(4.12)$$ Set
$${\cal
W}=[\sum_{i=1}^n(\mbb{C}x_i+\mbb{C}y_i)+\sum_{r=1}^m(\mbb{C}\sta_r+\mbb{C}\vt_r)]
[\sum_{j=1}^n(\mbb{C}\ptl_{x_{_j}}+\mbb{C}\ptl_{y_{_j}})+\sum_{s=1}^m(\mbb{C}\ptl_{\sta_s}+\mbb{C}\ptl_{\vt_s})].
\eqno(4.13)$$ Then
$$osp(2n|2m)|_{\cal A}=\{T\in{\cal W}\mid T(\eta)=0\}.\eqno(4.14)$$
Moreover,
$$\xi\Dlt=\Dlt \xi,\qquad \xi\eta=\eta
\xi\qquad\for\;\;\xi\in osp(2n|2m)\eqno(4.15)$$ as operators on
${\cal A}$. For $k\in\mbb{N}$, the subspace
$${\cal H}_k=\{f\in{\cal A}_k\mid\Dlt(f)=0\}\eqno(4.16)$$ forms an
$osp(2n|2m)$-submodule. First we prove the first conclusion in
Theorem 3:\psp

\pse

{\bf Theorem 4.1}. {\it Suppose $n>1$. For $k\in\mbb{N}$, ${\cal
H}_k$ is an irreducible $osp(2n|2m)$-module if and only if $k\leq
m+1-n$ or $k>2(m+1-n)$. When $k\leq m+1-n$, ${\cal
A}_k=\bigoplus_{i=0}^{\llbracket k/2 \rrbracket}\eta^i{\cal
H}_{k-2i}$ is a decomposition of irreducible
$osp(2n|2m)$-submodules}.

{\it Proof}. We take the subspace of diagonal matrices in
$osp(2n|2m)$ as a Carten subalgebras and the subspace
\begin{eqnarray*}& &osp(2n|2m)_+=\sum_{1\leq i<j\leq
n}[\mbb{C}(E_{i,j}-E_{n+j,n+i})+\mbb{C}(E_{i,n+j}-E_{j,n+i})] \\ & &
+\sum_{1\leq r<s\leq m}\mbb{C}
(E_{2n+r,2n+s}-E_{2n+m+s,2n+m+r})+\sum_{1\leq r\leq s\leq m}\mbb{C}
(E_{2n+r,2n+m+s}+E_{2n+s,2n+m+r})\\ & &+\sum_{i=1}^n\sum_{r=1}^m[
\mbb{C}(E_{i,2n+r}-E_{2n+m+r,n+i})+\mbb{C}(E_{i,2n+m+r}+E_{2n+r,n+i})]
\hspace{3cm}(4.17)\end{eqnarray*}
 as the space
spanned by positive root vectors. An $osp(2n|2m)$-{\it singular
vector} $v$ is a nonzero weight vector of $osp(2n|2m)$ such that
$osp(2n|2m)_+(v)=0$. We count singular vector up to a nonzero scalar
multiple.

 Observe ${\cal K}|_{\cal
A}=gl(n|m)|_{\cal A}$. According to the arguments in (2.51)-(2.60),
the homogeneous singular vectors of ${\cal K}$ are:
$$\{\eta^{\ell_3}x_1^{\ell_1}y_n^{\ell_2}\vec\vt_{s}\mid\ell_i\in\mbb{N};s\in\ol{1,m+1};
\ell_2(s-1)=0\}.\eqno(4.18)$$ By (4.10),
$(E_{n,2n+m+s}+E_{2n+s,2n})|_{\cal
A}=x_n\ptl_{\vt_s}+\sta_s\ptl_{y_n}$. Thus the homogeneous singular
vectors of $osp(2n|2m)$ are
$\{\eta^{\ell_2}x_1^{\ell_1}\mid\ell_1,\ell_2\in\mbb{N}\}.$
Moreover, (2.27) and (2.36) imply
$$\Dlt(\eta^{\ell_2}x_1^{\ell_1})=\ell_2(n-m+\ell_1+\ell_2-1)
=0\lra \ell_1+n\leq m\;\mbox{and}\;\ell_2=m+1-n-\ell_1.\eqno(4.19)$$
In this case, $\eta^{m+1-n-\ell_1}x_1^{\ell_1}\in{\cal
H}_{2(m+1-n)-\ell_1}$. Thus
$${\cal H}_k\;\mbox{has a unique singular vector if and only if}\;
k\leq m+1-n\;\mbox{or}\; k>2(m+1-n),\eqno(4.20)$$ and
$${\cal H}_k\;\mbox{has two singular vectors when}\;
 m+1-n< k\leq 2(m+1-n).\eqno(4.21)$$

Note $x_1^k\in{\cal H}_k$. Let $U$ be the $osp(2n|2m)$-submodule
generated by $x_1^k$. Repeatedly applying
$(E_{n+1,2n+r}-E_{2n+m+r,1})|_{\cal
A}=y_1\ptl_{\sta_r}-\vt_r\ptl_{x_1}$ (cf. (4.10)), we get
$$x_1^\ell\vec\vt_s\in U\qquad\for\;\;\ell+m+1-s=k.\eqno(4.22)$$
According to (4.8),  $(E_{n+1,n}-E_{2n,1})|_{\cal
A}=y_1\ptl_{x_n}-y_n\ptl_{x_1}$. Thus
$$x_1^{\ell_1}y_n^{\ell_2}\vec\vt_1\in
U\qquad\for\;\;\ell_1+\ell_2+m=k\eqno(4.23)$$ when $k\geq m$. Since
$${\cal H}_k=\sum_{i=0}^k{\cal H}_{i,k-i},\eqno(4.24)$$
Lemma 2.5, (4.22) and (4.23) imply $U={\cal H}_k$. So ${\cal H}_k$
is an $osp(2n|2m)$-module generated by $x_1^k$.

Suppose $k\leq m+1-n$ or $k>2(m+1-n)$. Let $M$ be a nonzero
$osp(2n|2m)$-submodule of ${\cal H}_k$. By (4.20), ${\cal H}_k$
contains a unique singular vector $x_1^k$. Thus $x_1^k\in M$. By the
above paragraph, ${\cal H}_k\subset M$. Hence ${\cal H}_k$ is
irreducible.

If ${\cal H}_k$ is irreducible, (4.21) implies $k\leq m+1-n$ or
$k>2(m+1-n)$.

Assume $k\leq m+1-n$. Note
$${\cal A}_j=\sum_{\ell=0}^j{\cal A}_{\ell,j-\ell}\qquad\for\;\;j\in\mbb{N}.\eqno(4.25)$$ By Theorem
1,
$${\cal A}_k=\bigoplus_{\ell=0}^k{\cal A}_{\ell,k-\ell}=\bigoplus_{\ell=0}^k\bigoplus_{i=0}^{\llbracket k/2 \rrbracket}\eta^i{\cal
H}_{\ell-i,k-\ell-i}=\bigoplus_{i=0}^{\llbracket k/2
\rrbracket}\eta^i{\cal H}_{k-2i}.\qquad\Box\eqno(4.26)$$
 \psp

 Fix $n_1,n_2\in\ol{1,n}$ with $n_1+1<n_2$.
  Changing operators $\ptl_{x_r}\mapsto -x_r,\;
 x_r\mapsto
\ptl_{x_r}$  for $r\in\ol{1,n_1}$ and $\ptl_{y_s}\mapsto -y_s,\;
 y_s\mapsto\ptl_{y_s}$  for $s\in\ol{n_2+1,n}$ in (4.6)-(4.10), we get a new representation of $osp(2n|2m)$ on
  ${\cal A}$ determined by
$$E_{2n+r,2n+s}|_{\cal A}=\sta_r\ptl_{\sta_s},\qquad
E_{2n+m+r,2n+m+s}|_{\cal A}=\vt_r\ptl_{\vt_s},\eqno(4.27)$$
$$E_{2n+r,2n+m+s}|_{\cal A}=\sta_r\ptl_{\vt_s},\qquad
E_{2n+m+r,2n+s}|_{\cal A}=\vt_r\ptl_{\sta_s},\eqno(4.28)$$
$$E_{i,j}|_{\cal A}=\left\{\begin{array}{ll}-x_j\ptl_{x_i}-\delta_{i,j}&\mbox{if}\;
i,j\in\ol{1,n_1};\\ \ptl_{x_i}\ptl_{x_j}&\mbox{if}\;i\in\ol{1,n_1},\;j\in\ol{n_1+1,n};\\
-x_ix_j &\mbox{if}\;i\in\ol{n_1+1,n},\;j\in\ol{1,n_1};\\
x_i\partial_{x_j}&\mbox{if}\;i,j\in\ol{n_1+1,n};
\end{array}\right.\eqno(4.29)$$
$$E_{n+i,n+j}|_{\cal A}=\left\{\begin{array}{ll}y_i\ptl_{y_j}&\mbox{if}\;
i,j\in\ol{1,n_2};\\ -y_iy_j&\mbox{if}\;i\in\ol{1,n_2},\;j\in\ol{n_2+1,n};\\
\ptl_{y_i}\ptl_{y_j} &\mbox{if}\;i\in\ol{n_2+1,n},\;j\in\ol{1,n_2};\\
-y_j\partial_{y_i}-\delta_{i,j}&\mbox{if}\;i,j\in\ol{n_2+1,n};
\end{array}\right.\eqno(4.30)$$
$$E_{i,n+j}|_{\cal
A}=\left\{\begin{array}{ll}
\ptl_{x_i}\ptl_{y_j}&\mbox{if}\;i\in\ol{1,n_1},\;j\in\ol{1,n_2};\\
-y_j\ptl_{x_i}&\mbox{if}\;i\in\ol{1,n_1},\;j\in\ol{n_2+1,n};\\
x_i\ptl_{y_j}&\mbox{if}\;i\in\ol{n_1+1,n},\;j\in\ol{1,n_2};\\
-x_iy_j&\mbox{if}\;i\in\ol{n_1+1,n},\;j\in\ol{n_2+1,n};\end{array}\right.\eqno(4.31)$$
$$E_{n+i,j}|_{\cal A}=\left\{\begin{array}{ll}
-x_jy_i&\mbox{if}\;j\in\ol{1,n_1},\;i\in\ol{1,n_2};\\
-x_j\ptl_{y_i}&\mbox{if}\;j\in\ol{1,n_1},\;i\in\ol{n_2+1,n};\\
y_i\ptl_{x_j}&\mbox{if}\;j\in\ol{n_1+1,n},\;i\in\ol{1,n_2};\\
\ptl_{x_j}\ptl_{y_i}&\mbox{if}\;j\in\ol{n_1+1,n},\;i\in\ol{n_2+1,n};
\end{array}\right.\eqno(4.32)$$
$$E_{i,2n+r}|_{\cal
A}=\left\{\begin{array}{ll}\ptl_{x_i}\ptl_{\sta_r}&\mbox{if}\;i\in\ol{1,n_1};\\
x_i\ptl_{\sta_r}&\mbox{if}\;i\in\ol{n_1+1,n};\end{array}\right.
E_{i,2n+m+r}|_{\cal
A}=\left\{\begin{array}{ll}\ptl_{x_i}\ptl_{\vt_r}&\mbox{if}\;i\in\ol{1,n_1};\\
x_i\ptl_{\vt_r}&\mbox{if}\;i\in\ol{n_1+1,n};\end{array}\right.\eqno(4.33)$$
$$E_{2n+r,i}|_{\cal
A}=\left\{\begin{array}{ll}-x_i\sta_r&\mbox{if}\;i\in\ol{1,n_1};\\\sta_r\ptl_{x_i}&\mbox{if}\;i\in\ol{n_1+1,n};\end{array}\right.
E_{2n+m+r,i}|_{\cal
A}=\left\{\begin{array}{ll}-x_i\vt_r&\mbox{if}\;i\in\ol{1,n_1};\\
\vt_r\ptl_{x_i}&\mbox{if}\;i\in\ol{n_1+1,n};\end{array}\right.\eqno(4.34)$$
$$E_{n+i,2n+r}|_{\cal
A}=\left\{\begin{array}{ll}y_i\ptl_{\sta_r}&\mbox{if}\;i\in\ol{1,n_2};\\
\ptl_{y_i}\ptl_{\sta_r}&\mbox{if}\;i\in\ol{n_2+1,n};\end{array}\right.\eqno(4.35)$$
$$E_{n+i,2n+m+r}|_{\cal
A}=\left\{\begin{array}{ll}y_i\ptl_{\vt_r}&\mbox{if}\;i\in\ol{1,n_2};\\
\ptl_{y_i}\ptl_{\vt_r}&\mbox{if}\;i\in\ol{n_2+1,n};\end{array}\right.\eqno(4.36)$$
$$E_{2n+r,n+i}|_{\cal
A}=\left\{\begin{array}{ll}\sta_r\ptl_{y_i}&\mbox{if}\;i\in\ol{1,n_2};\\
-y_i\sta_r&\mbox{if}\;i\in\ol{n_2+1,n};\end{array}\right.
\eqno(4.37)$$
$$E_{2n+m+r,n+i}|_{\cal
A}=\left\{\begin{array}{ll}\vt_r\ptl_{y_i}&\mbox{if}\;i\in\ol{1,n_2};\\
-y_i\vt_r&\mbox{if}\;i\in\ol{n_2+1,n};\end{array}\right.\eqno(4.38)$$
for $i,j\in\ol{1,n}$ and $r,s\in\ol{1,m}$.

The related Laplace operator becomes
$$\Dlt=-\sum_{i=1}^{n_1}x_i\ptl_{y_i}+\sum_{r=n_1+1}^{n_2}\ptl_{x_r}\ptl_{y_r}-\sum_{s=n_2+1}^n
y_s\ptl_{x_s}+\sum_{r=1}^m\ptl_{\sta_r}\ptl_{\vt_r}\eqno(4.39)$$ and
its dual
$$\eta=\sum_{i=1}^{n_1}y_i\ptl_{x_i}+\sum_{r=n_1+1}^{n_2}x_ry_r+\sum_{s=n_2+1}^n
x_s\ptl_{y_s}+\sum_{r=1}^m\sta_r\vt_r.\eqno(4.40)$$ It can be
verified that (4.15) holds again. Denote
\begin{eqnarray*}\qquad {\cal A}_{\la
k\ra}&=&\mbox{Span}\{x^\al
y^\be\Theta_1^{\ell_1'}\Theta_2^{\ell_2'}\mid\al,\be\in\mbb{N}^n;\ell_1',\ell_2'\in\mbb{N};\\
& &\sum_{r=n_1+1}^n\al_r-\sum_{i=1}^{n_1}\al_i+
\sum_{i=1}^{n_2}\be_i-\sum_{r=n_2+1}^n\be_r+\ell_1'+\ell_2'=\ell_2\}\hspace{2.6cm}(4.41)\end{eqnarray*}
for $k\in\mbb{Z}$.

Again we set ${\cal H}_{\la k\ra}=\{f\in {\cal A}_{\la k\ra}\mid
\Dlt(f)=0\}. $  Next we prove the second conclusion in Theorem
3:\psp

{\bf Theorem 4.2}. {\it Let $k\in\mbb{Z}$.  The $osp(2n|2m)$-module
${\cal H}_{\la k\ra}$ is irreducible if and only if $k\leq
n_1+m+1-n_2$. When $k\leq n_1+m+1-n_2$, ${\cal A}_{\la
k\ra}=\bigoplus_{i=0}^\infty\eta^i({\cal H}_{\la k-2i\ra})$ is the
decomposition of irreducible $osp(2n|2m)$-submodules.}

{\it Proof}. Observe ${\cal K}|_{\cal A}=gl(n|m)|_{\cal A}$ in terms
the representation of $gl(n|m)$ given in (1.18)-(1.22). Lemma 3.2
says that the homogeneous singular vectors of ${\cal K}$ are of the
form:
$$\eta^{\ell_1}(x_i^{m_1}y_j^{m_2}\vec\vt_s)\;\;
\mbox{with}\;\ell_1,m_1,m_2\in\mbb{N},\;s\in\ol{1,m+1}\eqno(4.42)$$
 and
$$(i,j)\in\{(n_1,n_2),
(n_1,n_2+1-\dlt_{n_2,n}),(n_1+1,n_2),(n_1+1,n_2+1-\dlt_{n_2,n})\}.\eqno(4.43)$$

{\it Claim 1}. For $k\in\mbb{N}$, ${\cal H}_{\la k\ra}$ is an
$osp(2n|2m)$-module generated by $x_{n_1+1}^k$ and ${\cal H}_{\la
-k\ra}$ is an $osp(2n|2m)$-module generated by $x_{n_1}^k$.\psp

Let $V$ be the $osp(2n|2m)$-module generated by $x_{n_1+1}^k\in
{\cal H}_{\la k\ra}$. By (4.31),
$$(E_{n_2+1,n+n_1+1}-E_{n_1+1,n+n_2+1})|_{\cal
A}=x_{n_2+1}\ptl_{y_{n_1+1}}+x_{n_1+1}y_{n_2+1}.\eqno(4.44)$$ Thus
$$(E_{n_2+1,n+n_1+1}-E_{n_1+1,n+n_2+1})^{k_1}(x_{n_1+1}^k)=x_{n_1+1}^{k+k_1}y_{n_2+1}^{k_1}\in
V.\eqno(4.45)$$ According to (4.32),
$$(E_{n+n_2,n_1+1}-E_{n+n_1+1,n_2})|_{\cal
A}=y_{n_2}\ptl_{x_{n_1+1}}-y_{n_1+1}\ptl_{x_{n_2}}.\eqno(4.46)$$
Repeatedly applying (4.46) to $x_{n_1+1}^k$, we obtain
$$x_{n_1+1}^{k_1}y_{n_2}^{k_2}\in
V\qquad\for\;k_1,k_2\in\mbb{N}\;\mbox{such
that}\;k_1+k_2=k.\eqno(4.47)$$ Note that (4.34) and (4.35) imply
$$(E_{n+n_1+1,2n+r}-E_{2n+m+r,n_1+1})|_{\cal
A}=y_{n_1+1}\ptl_{\sta_r}-\vt_r\ptl_{x_{n_1+1}}\eqno(4.48)$$
Applying (4.48) with various $r$ to (4.45) and (4.47), we obtain
$$x_{n_1+1}^{\ell_1}y_j^{\ell_2}\vec\vt_s\in {\cal H}_{\la
k\ra}\;\;\mbox{with}\;j\in\{n_2,n_2+1\}\lra
x_{n_1+1}^{\ell_1}y_j^{\ell_2}\vec\vt_s\in V.\eqno(4.49)$$ In terms
of (3.34),
$$v_{\ell,\ell'}\in V,\qquad \ell+\ell'=k.\eqno(4.50)$$
By Lemma 3.4,
$${\cal H}_{\la \ell,\ell'\ra}\subset V,\qquad \ell+\ell'=k.\eqno(4.51)$$
Therefore
$${\cal H}_{\la k\ra}=\bigoplus_{\ell,\ell'\in\mbb{Z};\;\ell+\ell'=k}
{\cal H}_{\la \ell,\ell'\ra}\subset V.\eqno(4.52)$$

Suppose $k>0$. Let $U$ be  the $osp(2n|2m)$-module generated by
$x_{n_1}^k{\cal H}_{\la -k\ra}$. Observe
$$(E_{n_2+1,n+n_1}-E_{n_1,n+n_2})|_{\cal
A}=x_{n_2+1}\ptl_{y_{n_1}}+y_{n_2+1}\ptl_{x_{n_1}}\eqno(4.53)$$ by
(4.31), and
$$(E_{n+n_1,n_2}-E_{n+n_2,n_1}|_{\cal
A}=y_{n_1}\ptl_{x_{n_1}}+x_{n_1}y_{n_1}\eqno(4.54)$$ by (4.32).
Repeatedly applying the above two equations to $x_{n_1}^k$, we have
$$x_{n_1}^{k_1}y_{n_2+1}^{k_2},x_{n_1}^{k+k_3}y_{n_2}^{k_3}\in
U\qquad\for\;\;k_1,k_2,k_3\in\mbb{N}\;\mbox{such
that}\;k_1+k_2=k.\eqno(4.55)$$ According to (4.34) and (4.35),
$$(E_{n+n_1,2n+r}-E_{2n+m+r,n_1})|_{\cal
A}=y_{n_1}\ptl_{\sta_r}+x_{n_1}\vt_r.\eqno(4.56)$$ Applying (4.56)
with various $r$ to (4.55), we find
$$x_{n_1}^{\ell_1}y_j^{\ell_2}\vec\vt_s\in {\cal H}_{\la
-k\ra}\;\;\mbox{with}\;j\in\{n_2,n_2+1\}\lra
x_{n_1}^{\ell_1}y_j^{\ell_2}\vec\vt_s\in V.\eqno(4.57)$$ In terms of
(3.34),
$$v_{\ell,\ell'}\in V,\qquad \ell+\ell'=-k.\eqno(4.58)$$
Lemma 3.4 gives
$${\cal H}_{\la \ell,\ell'\ra}\subset V,\qquad \ell+\ell'=-k.\eqno(4.59)$$
Thus
$${\cal H}_{\la -k\ra}=\bigoplus_{\ell,\ell'\in\mbb{Z};\;\ell+\ell'=-k}
{\cal H}_{\la \ell,\ell'\ra}\subset V.\eqno(4.60)$$ This prove Claim
1.\psp

{\it Claim 2}. For $n_1+m+1-n_2\geq k\in\mbb{Z}$,  any nonzero
$osp(2n|2m)$-submodule of ${\cal H}_{\la k\ra}$ contains
$x_{n_1+1}^k$ if $k\geq 0$ or $x_{n_1}^{-k}$ when $k<0$.\psp

Note
$$x_{n_1}^{m_1}y_{n_2}^{m_2}\vec\vt_s\in{\cal H}_{\la
m+m_2+1-s-m_1\ra},\qquad
x_{n_1}^{m_1}y_{n_2+1}^{m_2}\vec\vt_s\in{\cal H}_{\la
m+1-s-m_1-m_2\ra},\eqno(4.61)$$
$$x_{n_1+1}^{m_1}y_{n_2}^{m_2}\vec\vt_s\in{\cal H}_{\la
m+m_1+m_2+1-s\ra},\qquad
x_{n_1+1}^{m_1}y_{n_2+1}^{m_2}\vec\vt_s\in{\cal H}_{\la
m+m_1+1-s-m_2\ra}.\eqno(4.62)$$ For $\ell_1\in\mbb{N}+1$ and
$f\in{\cal H}_{\la k'\ra}$ with $k'\in\mbb{Z}$, (2.27) and (3.21)
imply
$$\Dlt\eta^{\ell_1}(f)=\ell_1(n_2-n_1-m+k'+\ell_1-1)\eta^{\ell_1-1}(f).\eqno(4.63)$$
Thus
$$\Dlt\eta^{\ell_1}(f)=0\dar k'\leq
n_1+m-n_2\;\;\mbox{and}\;\;\ell_1=n_1+m-n_2-k'+1.\eqno(4.64)$$ If
the condition holds, then
$$\eta^{\ell_1}(f)\in{\cal H}_{\la
2(n_1+m+1-n_2)-k'\ra}.\eqno(4.65)$$ Moreover,
$$2(n_1+m+1-n_2)-k'=n_1+m+2-n_2+(n_1+m-n_2-k')\geq n_1+m+2-n_2.\eqno(4.66)$$
This shows that
$${\cal H}_{\la k_1\ra}\bigcap(\bigcup_{i=1}^\infty \eta^i({\cal
H}))\neq \{0\}\dar k_1\geq  n_1+m+2-n_2.\eqno(4.67)$$

Suppose $k\leq n_1+m+1-n_2$. Then the singular vectors of ${\cal K}$
in ${\cal H}_{\la k\ra}$ are of the form
$$x_i^{m_1}y_j^{m_2}\vec\vt_s\;\;
\mbox{with}\;m_1,m_2\in\mbb{N},\;s\in\ol{1,m+1}\eqno(4.68)$$ with
(4.43). Observe that
$$(E_{n_1+1,n+n_2}-E_{n_2,n+n_1+1})|_{\cal
A}=x_{n_1+1}\ptl_{y_{n_2}}-x_{n_2}\ptl_{y_{n_1+1}}\eqno(4.69)$$ by
(4.29), and
$$(E_{n+n_2+1,n_1}-E_{n+n_1,n_2+1})|_{\cal
A}=-x_{n_1}\ptl_{y_{n_2+1}}-y_{n_1}\ptl_{x_{n_2+1}}\eqno(4.70)$$ by
(4.32). According to (4.33) and (4.37),
$$(E_{n_1+1,2n+m+r}+E_{2n+r,n+n_1+1})|_{\cal
A}=x_{n_1+1}\ptl_{\vt_r}+\sta_r\ptl_{y_{n_1+1}}.\eqno(4.71)$$ Let
$M$ be any nonzero $osp(2n|2m)$-submodule of ${\cal H}_{\la k\ra}$.
Then $M$ contains at least one of the ${\cal K}$-singular vectors in
(4.68). Applying (4.68)-(4.71), we get
$$x_{n_1}^{k_1}x_{n_2}^{k_2}\in M\;\;\mbox{for
some}\;k_1,k_2\in\mbb{N}\;\mbox{such that}\;k_2-k_1=k.\eqno(4.72)$$
By (4.29) and (4.30),
$$(E_{n_1,n_1+1}-E_{n+n_1+1,n+n_1})|_{\cal
A}=\ptl_{x_{n_1}}\ptl_{x_{n_1+1}}-y_{n_1+1}\ptl_{y_{n_1}}.\eqno(4.73)$$
Repeatedly applying (4.73) to (4.72), we obtain $x_{n_1+1}^k\in M$
if $k\geq 0$ or $x_{n_1}^{-k}\in M$ when $k<0$.\psp

The above claims show that ${\cal H}_{\la k\ra}$ is an irreducible
$osp(2n|2m)$-module if $k\leq  n_1+m+1-n_2$.

Note
$${\cal A}_{\la j\ra}=\sum_{\ell\in\mbb{Z}}{\cal A}_{\la\ell,j-\ell\ra},\qquad
{\cal H}_{\la j\ra}=\sum_{\ell\in\mbb{Z}}{\cal
H}_{\la\ell,j-\ell\ra}\qquad\for\;\;j\in\mbb{Z}.\eqno(4.74)$$ By
(4.67), $k\leq  n_1+m+1-n_2$ if ${\cal H}_{\la k\ra}$ is
irreducible. When $k\leq  n_1+m+1-n_2$,  Theorem 2 implies
$${\cal A}_{\la k\ra}=\bigoplus_{\ell\in\mbb{Z}}{\cal A}_{\la\ell,k-\ell\ra}=\bigoplus_{\ell\in\mbb{Z}}\bigoplus_{i=0}^\infty\eta^i{\cal
H}_{\la\ell-i,k-\ell-i\ra}=\bigoplus_{i=0}^\infty\eta^i{\cal H}_{\la
k-2i\ra}.\qquad\Box\eqno(4.75)$$
 \psp

We remark that ${\cal H}_{\la k\ra}$ is an indecomposable
$osp(2n|2m)$-module if $k \geq n_1+m+2-n_2$ by Claim 1 and (4.67).
This also implies that  ${\cal A}_{\la k\ra}$ is not completely
reducible $osp(2n|2m)$-module when $k\geq n_1+m+2-n_2$.

\section{Proof of Theorem 4}

\quad$\;\;$In this section, we want to prove Theorem 3.

Note
$$osp(2n+1|2m)_0=osp(2n|2m)_0+\sum_{i=1}^n[\mbb{C}(E_{0,i}-E_{n+i,0})
+\mbb{C}(E_{0,n+i}-E_{i,0})],\eqno(5.1)$$
\begin{eqnarray*}\qquad \qquad osp(2n+1|2m)_1&=&osp(2n|2m)_1+\sum_{r=1}^m[\mbb{C}(E_{0,2n+r}-E_{2n+m+r,0})
\\ & &+\mbb{C}(E_{0,2n+m+r}+E_{2n+r,0})].\hspace{4.9cm}(5.2)\end{eqnarray*} The Lie
superalgebra $osp(2n+1|2m)=osp(2n+1|2m)_0+osp(2n+1|2m)_1$ is a Lie
sub-superalgebra of $gl(2n+1|2m)$. Take settings in (1.28) and
(1.29). Now $osp(2n|2m)$ acts on ${\cal B}$ by the differential
operators in (4.6)-(4.10), namely, we change the subindex $|_{\cal
A}$ to $|_{\cal B}$. Extend the representation of $osp(2n|2m)$ on
${\cal B}$ to a representation of $osp(2n+1|2m)$ on ${\cal B}$ by:
$$(E_{0,i}-E_{n+i,0})|_{\cal B}=x_0\ptl_{x_i}-y_i\ptl_{x_0},
\qquad (E_{0,n+i}-E_{i,0})|_{\cal
B}=x_0\ptl_{y_i}-x_i\ptl_{x_0},\eqno(5.3)$$
$$(E_{0,2n+r}-E_{2n+m+r,0})|_{\cal B}=x_0\ptl_{\sta_r}-\vt_r\ptl_{x_0},
\;\; (E_{0,2n++m+r}+E_{2n+r,0})|_{\cal
B}=x_0\ptl_{\vt_r}+\sta_r\ptl_{x_0}\eqno(5.4)$$ for $i\in\ol{1,n}$
and $r\in\ol{1,m}$.

Set
$${\cal
W}'=[\sum_{i=0}^n(\mbb{C}x_i+\mbb{C}y_i)+\sum_{r=1}^m(\mbb{C}\sta_r+\mbb{C}\vt_r)]
[\sum_{j=0}^n(\mbb{C}\ptl_{x_{_j}}+\mbb{C}\ptl_{y_{_j}})+\sum_{s=1}^m(\mbb{C}\ptl_{\sta_s}+\mbb{C}\ptl_{\vt_s})].
\eqno(5.5)$$ Then
$$osp(2n+1|2m)|_{\cal B}=\{T\in {\cal W}'\mid
T(\eta')=0\}.\eqno(5.6)$$ Define
$${\cal H}'=\{f\in{\cal B}'\mid \Dlt'(f)=0\},\qquad {\cal
H}'_k={\cal H}'\bigcap {\cal B}_k.\eqno(5.7)$$

Again we take the subspace of diagonal matrices in $osp(2n+1|2m)$ as
a Cartan subalgebra and take the space spanned by positive roots:
$$osp(2n+1|2m)_+=\mbb{C}(E_{0,n+i}-E_{i,0})+\mbb{C}(E_{0,2n++m+r}+E_{2n+r,0})+
osp(2n|2m)_+.\eqno(5.8)$$ An $osp(2n+1|2m)$-{\it singular vector}
$v$ is a nonzero weight vector of $osp(2n+1|2m)$ such that
$osp(2n+1|2m)_+(v)=0$. We count singular vector up to a nonzero
scalar multiple. According to the proof of Theorem 4.1, any singular
vector of $osp(2n+1|2m)$ must be in $\mbb{C}[x_0,x_1,\eta]$, where
$$\eta=\sum_{i=1}^nx_iy_i+\sum_{r=1}^m\sta_r\vt_r.\eqno(5.9)$$
Note that $\eta'=x_0^2+2\eta$. By (5.8), $x_1^k$ is a singular
vector of $osp(2n+1|2m)$ for any $k\in\mbb{N}$. Thus a homogeneous
singular vector of $osp(2n+1|2m)$ must be of the form
$$f=\sum_{i=0}^\ell
b_ix_0^{2i+\iota}\eta^{\ell-i}x_1^k,\eqno(5.10)$$ where
$b_i\in\mbb{C}$, $\ell,k\in\mbb{N}$ and $\iota=0,1.$ Note
$$(E_{0,n+i}-E_{i,0})(f)=(x_0\ptl_{y_i}-x_i\ptl_{x_0})(f)=0\dar
f=b_0(\eta')^\ell x_1^k.\eqno(5.11)$$ Thus $\{(\eta')^\ell x_1^k\mid
\ell,k\in\mbb{N}\}$ are all the homogeneous singular vectors of
$osp(2n+1|2m)$ in ${\cal B}$.

Observe that
$$[\Dlt',\eta']=2+4(n-m)+4[x_0\ptl_{x_0}+\sum_{i=1}^n(x_i\ptl_{x_i}
+y_i\ptl_{y_i})+\sum_{r=1}^m(\sta_r\ptl_{\sta_r}+\vt_r\ptl_{\vt_r})]\eqno(5.12)$$
 by (2.27) and (2.36). So
$$\Dlt'((\eta')^\ell g)=2\ell(1+2(n-m+k+\ell-1))(\eta')^{\ell-1}
g\qquad\for\;\;g\in{\cal H}_k\eqno(5.13)$$ Thus
$${\cal H}'_k\;\mbox{has a unique singular vector}\;x_1^k\;\mbox{for
any}\;k\in\mbb{N}.\eqno(5.14)$$ Indeed we have the first conclusion
in Theorem 4:\psp

{\bf Theorem 5.1}. {\it For any $k\in\mbb{N}$, ${\cal H}_k'$ is an
irreducible $osp(2n+1|2m)$-module. Moreover, ${\cal
B}=\bigoplus_{\ell,k=0}^\infty (\eta')^\ell {\cal H}_k$ is a direct
sum of irreducible $osp(2n+1|2m)$-submodules.}

{\it Proof}. By the arguments in (3.52)-(3.58), we only need to
prove that ${\cal H}'_k$ is an $osp(2n+1|2m)$-module generated by
$x_1^k$. For $\iota=0,1$, we define
$$T_\iota=\sum_{i=0}^\infty\frac{(-2)^ix_0^{2i+\iota}}{(2i+\iota)!}\Dlt^i,\qquad\Dlt=\sum_{i=1}^n\ptl_{x_i}\ptl_{y_i}
+\sum_{r=1}^m\ptl_{\sta_r}\ptl_{\vt_r} .\eqno(5.15)$$ We take the
notations in (2.104) and (2.105). By Lemma 2.6,
\begin{eqnarray*}\qquad{\cal H}_k'&=&\mbox{Span}\{T_\iota(x^\al y^\be
\sta_{\vec j}\vt_{\vec j'})\mid\al,\be\in\mbb{N}^n;\vec j\in
\G_{k_1},\;\vec j'\in
\G_{k_2};\\
&&\qquad\iota\in\{0,1\};|\al|+|\be|+k_1+k_2+\iota=k\}.\hspace{4.9cm}(5.16)\end{eqnarray*}
Let $U$ be the $osp(2n+1|2m)$-module generated by $x_1^k\in{\cal
H}_k$. Since $o(2n+1,\mbb{C})$ is a subalgebra of $osp(2n+1|2m)$,
the known results of the representation of $o(2n+1,\mbb{C})$ on
$\mbb{C}[x_0,x_1,...,x_{2n}]$ show
$$T_\iota(x^\al y^\be)\in
U\qquad\for\;\;\al,\be\in\mbb{N}^n;|\al|+|\be|=k.\eqno(5.17)$$
Repeatedly applying $(E_{2n+r,i}+E_{n+i,2n+m+r})|_{\cal
A}=\sta_r\ptl_{x_i}+y_i\ptl_{\vt_r}$
 to (5.17) with $i\in\ol{1,n}$
and $r\in\ol{1,m}$, we obtain
$$T_\iota(x^\al y^\be\sta_{\vec j})\in
U\qquad\for\;\;\al,\be\in\mbb{N}^n,\;\vec
j\in\G_{k_1};|\al|+|\be|+k_1=k.\eqno(5.18)$$ Finally, we get
$U={\cal H}_k'$ by repeatedly applying
$(E_{i,2n+r}-E_{2n+m+r,n+i})|_{\cal
B}=x_i\ptl_{\sta_r}-\vt_r\ptl_{y_i}$ to (5.18) with $i\in\ol{1,n}$
and $r\in\ol{1,m}$.$\qquad\Box$\psp

 Next
$osp(2n|2m)$ acts on ${\cal B}$ via the differential operators in
(4.27)-(4.38), namely, we change the subindex $|_{\cal A}$ to
$|_{\cal B}$. Moreover, we extend the representation of $osp(2n|2m)$
on ${\cal B}$ to a representation of $osp(2n+1|2m)$ on ${\cal B}$
by:
$$(E_{0,i}-E_{n+i,0})|_{\cal
B}=\left\{\begin{array}{ll}-x_0x_i-y_i\ptl_{x_0}&\mbox{if}\;i\in\ol{1,n_1},\\
x_0\ptl_{x_i}-y_i\ptl_{x_0}&\mbox{if}\;i\in\ol{n_1+1,n_2},\\
x_0\ptl_{x_i}-\ptl_{x_0}\ptl_{y_i}&\mbox{if}\;i\in\ol{n_2+1,n};\end{array}\right.
\eqno(5.19)$$
$$(E_{0,n+i}-E_{n,0})|_{\cal
B}=\left\{\begin{array}{ll}x_0\ptl_{y_i}-\ptl_{x_i}\ptl_{x_0}&\mbox{if}\;i\in\ol{1,n_1},\\
x_0\ptl_{y_i}-x_i\ptl_{x_0}&\mbox{if}\;i\in\ol{n_1+1,n_2},\\
-x_0y_i-x_i\ptl_{x_0}&\mbox{if}\;i\in\ol{n_2+1,n};\end{array}\right.
\eqno(5.20)$$
$$(E_{0,2n+r}-E_{2n+m+r,0})|_{\cal B}=x_0\ptl_{\sta_r}-\vt_r\ptl_{x_0},
\;\; (E_{0,2n++m+r}+E_{2n+r,0})|_{\cal
B}=x_0\ptl_{\vt_r}+\sta_r\ptl_{x_0}\eqno(5.21)$$ for $i\in\ol{1,n}$
and $r\in\ol{1,m}$.

Now the corresponding Laplace operator becomes
$$\Dlt'=\ptl_{x_0}^2+2\Dlt,\;\;\Dlt=-\sum_{i=1}^{n_1}x_i\ptl_{y_i}+\sum_{r=n_1+1}^{n_2}\ptl_{x_r}\ptl_{y_r}-\sum_{s=n_2+1}^n
y_s\ptl_{x_s}+\sum_{r=1}^m\ptl_{\sta_r}\ptl_{\vt_r}\eqno(5.22)$$ and
its dual
$$\eta'=x_0^2+2\eta,\;\;\eta=\sum_{i=1}^{n_1}y_i\ptl_{x_i}+\sum_{r=n_1+1}^{n_2}x_ry_r+\sum_{s=n_2+1}^n
x_s\ptl_{y_s}+\sum_{r=1}^m\sta_r\vt_r.\eqno(5.23)$$ We take the
notation in (4.41) and set
$${\cal B}_{\la k\ra}=\sum_{i=0}^\infty {\cal A}_{\la
k-i\ra}x_0^i,\qquad {\cal H}'_{\la k\ra}=\{f\in {\cal B}_{\la
k\ra}\mid \Dlt'(f)=0\}.\eqno(5.24)$$ Then we have the second
conclusion in Theorem 4:\psp

{\bf Theorem 5.2}.  {\it For any $k\in\mbb{Z}$, ${\cal H}_{\la
k\ra}'$ is an irreducible $osp(2n+1|2m)$-module. Moreover, ${\cal
B}=\bigoplus_{\ell,k=0}^\infty (\eta')^\ell {\cal H}_{\la k\ra}$ is
a direct sum of irreducible $osp(2n+1|2m)$-submodules.}

{\it Proof}. We define $T_\iota$ as in (5.15) with $\Dlt$ in (5.22).
By Lemma 2.6,
$${\cal H}'_{\la k\ra}=T_0({\cal A}_{\la k\ra})+T_1({\cal A}_{\la
k-1\ra})\qquad\for\;\;k\in\mbb{Z}.\eqno(5.25)$$ Since
$\Dlt\xi=\xi\Dlt$ for $\xi\in osp(2n|2m)$, we have
$$\xi(T_\iota(f))=T_\iota(\xi(f))\qquad\for\;\;\xi\in osp(2n|2m),\;
f\in{\cal A}.\eqno(5.26)$$

First we consider ${\cal H}'_{\la k\ra}$ with $k\in\mbb{N}$. Let $V$
be any nonzero $osp(2n+1|2m)$-submodule of ${\cal H}'_{\la k\ra}$.
According to the arguments in paragraph of (4.68)-(4.73), $V$
contains some $T_\iota(\eta^\ell(x_{n_1+1}^{k-\iota-2\ell})).$
According to (5.20),
$$(E_{n_1+1,0}-E_{0,n+n_1+1})|_{\cal B}=x_{n_1+1}\ptl_{x_0}-x_0\ptl_{y_{n_1+1}}.\eqno(5.27)$$
Moreover, as operators on ${\cal B}$,
\begin{eqnarray*}& &[E_{n_1+1,0}-E_{0,n+n_1+1},T_0]\\&=&
[x_{n_1+1}\ptl_{x_0}-x_0\ptl_{y_{n_1+1}},\sum_{i=0}^\infty\frac{(-2)^ix_0^{2i}}{(2i)!}\Dlt^i]\\
&=&[x_{n_1+1},\sum_{i=1}^\infty\frac{(-2)^ix_0^{2i}}{(2i)!}\Dlt^i]\ptl_{x_0}+x_{n_1+1}
\sum_{i=1}^\infty\frac{(-2)^ix_0^{2i-1}}{(2i-1)!}\Dlt^i\\&=&
-[\sum_{i=1}^\infty\frac{i(-2)^ix_0^{2i}}{(2i)!}\Dlt^{i-1}\ptl_{x_0}
+\sum_{i=1}^\infty\frac{i(-2)^ix_0^{2i-1}}{(2i-1)!}\Dlt^{i-1}]\ptl_{y_{n_1+1}}
-2(T_1\Dlt) x_{n_1+1},\hspace{1.5cm}(5.28)
\end{eqnarray*}
\begin{eqnarray*}& &[E_{n_1+1,0}-E_{0,n+n_1+1},T_1]\\&=&
[x_{n_1+1}\ptl_{x_0}-x_0\ptl_{y_{n_1+1}},\sum_{i=0}^\infty\frac{(-2)^ix_0^{2i+1}}{(2i+1)!}\Dlt^i]\\
&=&[x_{n_1+1},\sum_{i=1}^\infty\frac{(-2)^ix_0^{2i+1}}{(2i+1)!}\Dlt^i]\ptl_{x_0}+x_{n_1+1}
\sum_{i=0}^\infty\frac{(-2)^ix_0^{2i}}{(2i)!}\Dlt^i\\&=&
-[\sum_{i=1}^\infty\frac{i(-2)^ix_0^{2i+1}}{(2i+1)!}\Dlt^{i-1}\ptl_{x_0}
+\sum_{i=1}^\infty\frac{i(-2)^ix_0^{2i}}{(2i)!}\Dlt^{i-1}]\ptl_{y_{n_1+1}}
+T_0 x_{n_1+1}.\hspace{2.4cm}(5.29)
\end{eqnarray*}
If $T_0(\eta^\ell(x_{n_1+1}^{k-2\ell}))\in V$ for some
$\ell\in\mbb{N}+1$, we have
\begin{eqnarray*}& &(E_{n_1+1,0}-E_{0,n+n_1+1})T_0(\eta^\ell(x_{n_1+1}^{k-2\ell}))\\
&=&
([E_{n_1+1,0}-E_{0,n+n_1+1},T_0]+T_0(E_{n_1+1,0}-E_{0,n+n_1+1}))(\eta^\ell(x_{n_1+1}^{k-2\ell}))
\\&=&-[\sum_{i=1}^\infty\frac{i(-2)^ix_0^{2i-1}}{(2i-1)!}\Dlt^{i-1}\ptl_{y_{n_1+1}}+x_0T_0\ptl_{y_{n_1+1}}
+2(T_1\Dlt)x_{n_1+1}](\eta^\ell(x_{n_1+1}^{k-2\ell}))
\\&=& [T_1\ptl_{y_{n_1+1}}-2(T_1\Dlt)x_{n_1+1}](\eta^\ell(x_{n_1+1}^{k-2\ell}))
\\ &=&\ell[1-2(m+n_1-n_2+\ell-1)](\eta^{\ell-1}(x_{n_1+1}^{k-2(\ell-1)-1}))
\\ &=&\ell[3-2(m+n_1-n_2+\ell)]T_1(\eta^{\ell-1}(x_{n_1+1}^{k-2(\ell-1)-1}))\in
V\hspace{4.9cm}(5.30)\end{eqnarray*} by (4.63), (5.27) and (5.28).
So $T_1(\eta^{\ell-1}(x_{n_1+1}^{k-2(\ell-1)-1}))\in V$. When
$T_1(\eta^\ell(x_{n_1+1}^{k-2\ell-1}))\in V$ for some
$\ell\in\mbb{N}$, (5.27) and (5.29) yield
$$(E_{n_1+1,0}-E_{0,n+n_1+1})T_1(\eta^\ell(x_{n_1+1}^{k-2\ell-1}))=T_0(\eta^\ell(x_{n_1+1}^{k-2\ell}))\in
V.\eqno(5.31)$$ By induction on $\ell$, we have
$x_{n_1+1}^k=T_0(x_{n_1+1}^k)\in V$.

Note
$$(E_{n+i,n+n_1+1}-E_{n_1+1,i})|_{\cal
B}=y_i\ptl_{y_{n_1+1}}+x_ix_{n_1+1}\qquad\for\;\;i\in\ol{1,n_1}\eqno(5.32)$$
and
$$(E_{n_2+r,n_2}-E_{n+n_2,n+m_2+r})|_{\cal
A}=x_{n_2+r}\ptl_{x_{n_2}}+y_{n_2}y_{n_2+r}\;\;\for\;\;r\in\ol{1,n-n_2}\;\;\mbox{if}\;\;n_2<n\eqno(5.33)$$
by (4.29) and (4.30). Repeatedly applying (5.32) and (5.33) to
(5.31) with various $i\in\ol{1,n_1}$ and $r\in\ol{1,n-n_2}$ if
$n_2<n$, we have
$$[\prod_{i=1}^{n_1+1}x_i^{\al_i}][\prod_{j=n_2}^ny_j^{\be_j}]\in
V\qquad\for\;\;\al_i,\be_j\in\mbb{N};\al_{n_1+1}+\be_{n_2}-\sum_{i=1}^{n_1}\al_i-\sum_{r=n_2+1}^n\be_r=k.\eqno(5.34)$$

Denote
$$I=\{0,\ol{n_1+1,n_2},\ol{n+n_1+1,n+n_2},\ol{2n+1,2n+2m}\}.\eqno(5.35)$$
Then the Lie subalgebra
$${\cal G}=osp(2n+1|2m)\bigcap(\sum_{i,j\in I}\mbb{C}E_{i,j})\cong
osp(2(n_2-n_1)+1|2m).\eqno(5.36)$$ Applying Theorem 5.1 to ${\cal
G}$ and
$\mbb{C}[x_0,x_{n_1+1},...,x_{n_2},y_{n_1+1},...,y_{n_2},\sta_1,...,\sta_m,\vt_1,...,\vt_m]$,
we get
$$T_\iota(x^\al y^\be\sta_{\vec j}\vt_{\vec j'})\in V\eqno(5.37)$$
for $\al,\be\in\mbb{N}^n$, $\vec j\in\G_{k_1}$ and $\vec
j'\in\G_{k_2}$ such that $\be_i=0$ if $i\leq n_1$ and $\al_j=0$ if
$j>n_2$, and
$$\iota+k_1+k_2+\sum_{r=n_1+1}^{n_2}(\al_r+\be_r)-\sum_{i=1}^{n_1}\al_i-\sum_{j=n_2+1}^n\be_j=k.\eqno(5.38)$$

Repeatedly applying (5.32) to (5.38) under above conditions with
various $i\in\ol{1,n_1}$, we obtain (5.38) for
$\al,\be\in\mbb{N}^n$, $\vec j\in\G_{k_1}$ and $\vec j'\in\G_{k_2}$
such that $\al_i=0$ if $i>n_2$, and
$$\iota+k_1+k_2+\sum_{r=n_1+1}^{n_2}\al_r+\sum_{s=1}^{n_2}\be_s-\sum_{i=1}^{n_1}\al_i-\sum_{j=n_2+1}^n\be_j=k.\eqno(5.39)$$
 Observe
$$(E_{n_2+r,n_1+s}-E_{n+n_1+s,n+n_2+r})|_{\cal
B}=y_{n_1+s}y_{n_2+r}+x_{n_2+r}\ptl_{x_{n_1+s}}\eqno(5.40)$$ for
$r\in\ol{1,n-n_2}$ and $s\in\ol{1,n_2-n_1}$ by (4.29) and (4.30).
Repeatedly applying (5.40) to (5.38) with $\al_i=0$ if $i>n_2$, we
obtain ${\cal H}_{\la k\ra}\subset V$ by (5.25). So ${\cal H}_{\la
k\ra}$ is an irreducible $osp(2n+1|2m)$-module.

Next we consider ${\cal H}_{\la -k\ra}$ with $k\in\mbb{N}+1$.  Let
$U$ be any nonzero $osp(2n+1|2m)$-submodule of ${\cal H}'_{\la
-k\ra}$. According to the arguments in paragraph of (4.68)-(4.73),
$U$ contains some $T_\iota(\eta^\ell(x_{n_1}^{k+\iota+2\ell})).$
Observe Note
$$(E_{n_1,0}-E_{0,n+n_1})=\ptl_{x_0}\ptl_{x_{n_1}}-x_0\ptl_{y_{n_1}}\eqno(5.41)$$
by (4.29) and (4.30). As operators on ${\cal B}$,
\begin{eqnarray*}& &[E_{n_1,0}-E_{0,n+n_1},T_0]\\&=&
[\ptl_{x_0}\ptl_{x_{n_1}}-x_0\ptl_{y_{n_1}}
,\sum_{i=0}^\infty\frac{(-2)^ix_0^{2i}}{(2i)!}\Dlt^i]\\
&=&\sum_{i=1}^\infty\frac{(-2)^ix_0^{2i-1}}{(2i-1)!}\Dlt^i\ptl_{x_{n_1}}
-\ptl_{x_0}\sum_{i=1}^\infty\frac{i(-2)^ix_0^{2i}}{(2i)!}\Dlt^{i-1}\ptl_{y_{n_1}}
\\&=&-2T_1\Dlt\ptl_{x_{n_1}}-
\sum_{i=1}^\infty\frac{i(-2)^ix_0^{2i-1}}{(2i-1)!}\Dlt^{i-1}\ptl_{y_{n_1}}-
\sum_{i=1}^\infty\frac{i(-2)^ix_0^{2i}}{(2i)!}\Dlt^{i-1}\ptl_{y_{n_1}}\ptl_{x_0},\hspace{1.7cm}(5.42)
\end{eqnarray*}
\begin{eqnarray*}& &[E_{n_1,0}-E_{0,n+n_1},T_1]\\&=&
[\ptl_{x_0}\ptl_{x_{n_1}}-x_0\ptl_{y_{n_1}},\sum_{i=0}^\infty\frac{(-2)^ix_0^{2i+1}}{(2i+1)!}\Dlt^i]\\
&=&\sum_{i=0}^\infty\frac{(-2)^ix_0^{2i}}{(2i)!}\Dlt^i\ptl_{x_{n_1}}-
\ptl_{x_0}\sum_{i=1}^\infty\frac{i(-2)^ix_0^{2i+1}}{(2i+1)!}\Dlt^{i-1}\ptl_{y_{n_1}}
\\ &=& T_0\ptl_{x_{n_1}}-\sum_{i=1}^\infty\frac{i(-2)^ix_0^{2i}}{(2i)!}\Dlt^{i-1}\ptl_{y_{n_1}}
-\sum_{i=1}^\infty\frac{i(-2)^ix_0^{2i+1}}{(2i+1)!}\Dlt^{i-1}\ptl_{y_{n_1}}\ptl_{x_0}
.\hspace{2.7cm}(5.43)
\end{eqnarray*}

If $T_0(\eta^\ell(x_{n_1}^{k+2\ell}))\in U$ for some
$\ell\in\mbb{N}+1$, we have
\begin{eqnarray*}& &(E_{n_1,0}-E_{0,n+n_1})T_0(\eta^\ell(x_{n_1+1}^{k-2\ell}))
\\&=& [T_1\ptl_{y_{n_1}}-2T_1\Dlt\ptl_{x_{n_1}}](\eta^\ell(x_{n_1+1}^{k+2\ell}))
\\ &=&(k+2\ell)\ell[1+2(n_1+k+\ell-n_2)]T_1(\eta^{\ell-1}(x_{n_1+1}^{k+2\ell-1}))
\in V\hspace{4.1cm}(5.44)\end{eqnarray*} by (4.63), (5.41) and
(5.42). So $T_1(\eta^{\ell-1}(x_{n_1}^{k+2(\ell-1)+1}))\in U$. When
$T_1(\eta^\ell(x_{n_1}^{k+2\ell+1}))\in V$ for some
$\ell\in\mbb{N}$, (5.41) and (5.43) yield
$$(E_{n_1,0}-E_{0,n+n_1})T_1(\eta^\ell(x_{n_1}^{k+2\ell+1}))=(k+2\ell+1)T_0(\eta^\ell(x_{n_1}^{k+2\ell}))\in
U.\eqno(5.45)$$ By induction on $\ell$, we have
$x_{n_1}^k=T_0(x_{n_1}^k)\in U$.

According to (5.32) with $i=n_1$,
$$x_{n_1}^{k+k'}x_{n_1+1}^{k'}\in
U\qquad\for\;\;k'\in\mbb{N}.\eqno(5.46)$$
 Moreover,
$$(E_{i,n_1}-E_{n+n_1,n+i})|_{\cal
B}=x_i\ptl_{x_{n_1}}-y_{n_1}\ptl_{y_i}\qquad\for
i\in\ol{1,n_1-1}\eqno(5.47)$$ by (4.29) and (4.30). Repeatedly
applying (5.47) to (5.46) with various $i\in\ol{1,n_1-1}$, we have
$$\prod_{i=1}^{n_1+1}x_i^{\al_i}\in
U\qquad\for\;\;\al_i\in\mbb{N};\al_{n_1+1}-\sum_{i=1}^{n_1}\al_i=-k.\eqno(5.48)$$
Observe
$$(E_{n_2+r,n+1}-E_{1,n+n_2+r})|_{\cal
A}=x_{n_2+r}\ptl_{y_1}+y_{n_2+r}\ptl_{x_1}\;\;\for\;\;r\in\in\ol{1,n-n_2}\;\;\mbox{if}\;\;n_2<n\eqno(5.49)$$
by (4.31). Repeatedly applying (5.49) to (5.48) with various
$r\in\ol{1,n-n_2}$ if $n_2<n$, we find
$$[\prod_{i=1}^{n_1+1}x_i^{\al_i}][\prod_{j=n_2+1}^ny_j^{\be_j}]\in
U\qquad\for\;\;\al_i,\be_j\in\mbb{N};\al_{n_1+1}-\sum_{i=1}^{n_1}\al_i-\sum_{j=n_2+1}^n=-k.\eqno(5.50)$$
By the same arguments from (5.36) to the end of the paragraph below
(5.40) with $k$ replaced by $-k$, we prove that ${\cal H}_{\la
-k\ra}$ is an irreducible $osp(2n+1|2m)$-module.

We calculate
$$[\Dlt',\eta']=2+4(n_2-n_1-m)+4[x_0\ptl_{x_0}+\sum_{i=1}^n(x_i\ptl_{x_i}
+y_i\ptl_{y_i})+\sum_{r=1}^m(\sta_r\ptl_{\sta_r}+\vt_r\ptl_{\vt_r})].\eqno(5.51)$$
By the arguments in (3.52)-(3.58),  ${\cal
B}=\bigoplus_{\ell,k=0}^\infty (\eta')^\ell {\cal H}_{\la k\ra}$ is
a direct sum of irreducible $osp(2n+1|2m)$-submodules for any
$k\in\mbb{Z}.\qquad\Box$\psp

\vspace{1cm}

\noindent{\Large \bf References}

\hspace{0.5cm}

\begin{description}

\item[{[C]}] B. Cao, Solutions of Navier Equations and Their
Representation Structure, {\it Adv. Appl. Math.} \textbf{43} (2009),
331-374.

\item[{[DES]}] M. Davidson, T.
Enright, and R. Stanke, {\it Differential Operators and Highest
Weight Representations}, Memoirs of American Mathematical Society
{\bf 94}, no. 455, 1991.

\item[{[FC]}] F. M. Fern\'{a}ndez and E. A. Castro, {\it Algebraic
Methods in Quantum Chemistry and Physics}, CRC Press, Inc., 1996.

\item[{[FSS]}] L.Frappat, A.Sciarrino and P.sorba, Dictionary on Lie
Algebras and Superalgebras, Academic Press,2000

\item[{[G]}] H. Georgi, {\it Lie Algebras in Particle Physics},
Second Edition, Perseus Books Group, 1999.

\item[{[Ho]}] R. Howe, Perspectives on invariant theory: Schur
duality, multiplicity-free actions and beyond, {\it The Schur
lectures} (1992) ({\it Tel Aviv}), 1-182, {\it Israel Math. Conf.
Proc.,} 8, {\it Bar-Ilan Univ., Ramat Gan,} 1995.

\item[{[Hu]}] J. E. Humphreys, {\it Introduction to Lie Algebras and Representation Theory},
 Springer-Verlag New York Inc., 1972.

\item[{[L]}] F. S. Levin, {\it An Introduction to Quantum Theory},
Cambridge University Press, 2002.

\item[{[LX]}] C. Luo and X. Xu,  Oscillator variations of the classical theorem on harmonic polynomials, {\it arXiv:1012.2391v1[math.RT].}

 \item[{[LF]}] W. Ludwig and C. Falter, {\it
Symmetries in Physics}, Second Edition, Springer-Verlag,
Berlin/Heidelberg, 1996.

 \item[{[O]}] P. Olver, Conservation laws in elasticity, II, Linear
homogeneous isotropic elastostatics, \emph{Arch. Rational Mech.
Anal.} \textbf{85} (1984), no.2, 131-160.

\item[{[X1]}] X.Xu, Flag partial differential equations and
representations of Lie algebras, {\it Acta Appl Math} {\bf 102}
(2008), 249--280.

\item[{[X2]}] X. Xu, A cubic $E_6$-generalization of the classical theorem on harmonic
polynomials, {\it J. Lie Theory} {\bf 21} (2011), 145-164.

\item[{[Z]}] R. B. Zhang, Orthosymplectic Lie superalgebras in
superspace analogues of quantum Kepler problems, {\it Commun. Math.
Phys.} {\bf 280} (2008), 545-562.

\end{description}

\end{document}